\documentclass[11pt]{article}

\usepackage{amsmath}
\usepackage{amssymb}
\usepackage{theorem}
\usepackage{euscript}
\usepackage{amscd}

\usepackage{pstricks}
\usepackage{pst-node,pst-tree}
\usepackage[arrow,curve,matrix,tips,frame]{xy}
\newpsobject{showgrid}{psgrid}{subgriddiv=1,griddots=10,gridlabels=6pt}

\usepackage{ifthen}

\renewcommand{\title}[1]{
  \addvspace{3\baselineskip}  
  \begin{center} \LARGE \bf #1
  \end{center}
  \addvspace{2\baselineskip}}   

\renewcommand{\author}[1]{
  \addvspace{-1\baselineskip}  
  \begin{center} \large \sc #1
  \end{center}
  \addvspace{2\baselineskip}}   

\makeatletter

\def\section{%
        \@startsection{section}{1}{\z@}%
        {8ex plus 6ex minus 3ex}{\baselineskip}%
        {\normalfont\large\scshape\centering}%
        }

\renewcommand{\paragraph}[1]{{\par\removelastskip\vskip.5\baselineskip%
         \indent{\scshape{#1}}{\ifperiod.\else\global\periodtrue\fi}%
         \rm \ignorespaces}}

\let\goth=\mathfrak
\let\calligraphy=\mathcal

\def\CC{{\mathbb C}}

\def\PP{{\mathbb P}}
\def\QQ{{\mathbb Q}}

\def\ZZ{{\mathbb Z}}

\def\Ff{{\calligraphy F}}
\def\Gg{{\calligraphy G}}

\def\Jj{{\calligraphy J}}

\def\Nn{{\calligraphy N}}
\def\Oo{{\calligraphy O}}

\def\Tt{{\calligraphy T}}

\def\MMM{{\goth M}}

\def\RRR{{\goth R}}

\def\hbar{{\,\overline{\!h}}}

\def\Ctilde{{\,\widetilde{\!C}}}

\def\sS{{\boldsymbol{S}}}

\def\card{\operatorname{card}}

\def\coker{\operatorname{coker}}

\def\Ext{\operatorname{Ext}}

\def\Gr{\operatorname{Gr}}

\let\lra=\longrightarrow

\def\ie{\idestdefinition}
\def\eg{{\it e.g.}~}
\def\cf{{\it cf.}~}

\def\inv{^{-1}}

\let\phi=\varphi
\let\epsilon=\varepsilon

\newcommand{\textq}[1]{\quad\text{#1}\quad}

\newcommand{\listspace}{\setlength{\itemsep}{-.5pt}}

\newcounter{icounter}

\def\theoname{Theorem}
\def\lemmaname{Lemma}
\def\propositionname{Proposition}
\def\notationname{Notation}
\def\corollaryname{Corollary}
\def\conjecturename{Conjecture}
\def\remarkname{Remark}
\def\remarksname{Remarks}
\def\examplename{Example}
\def\examplesname{Examples}
\def\definitionname{Definition}
\def\definitionsname{Definitions}
\def\notationname{Notation}

\def\proofname{Proof}
\def\exercisename{Exercise}

\def\idestdefinition{{\it i.e.}~}

\def\Dquad{\hskip 0.6em plus .02em minus .2em}  
\def\Dpar{\belowdisplayskip=0pt\belowdisplayshortskip=0pt\par}

\def\bigpenalty{\interlinepenalty=\@M}
\def\smallpenalty{\interlinepenalty=100}

\newif\ifperiod \periodtrue 

\def\D@makemargins{%
  \labelsep=0pt
  \itemindent=0pt
  \labelwidth=0pt}

\def\D@restoremargins{%
  \labelsep=5pt
  \itemindent=0pt
  \leftmargin=5mm  
  \labelwidth=\leftmargin \advance\labelwidth by -\labelsep}

\def\th@Dindent{\hspace\parindent}
\def\th@Dheadingshape{\scshape}

\gdef\th@DthAndSuchtheo{%
  \D@makemargins%
  \def\@begintheorem##1##2{%
  \item[]\th@Dindent{\th@Dheadingshape ##1~\rm ##2.}\Dquad         
        \D@restoremargins}%
  \def\@opargbegintheorem##1##2##3{\def\next{##3}%
  \item[]\th@Dindent{\th@Dheadingshape ##1~\rm ##2\ifx\next\empty
  \else\ {\normalfont(##3)}\fi.}         
        \D@restoremargins}}

\gdef\th@DthAndSuchtheostar{%
  \D@makemargins%
  \def\@begintheorem##1##2{%
  \item[]\th@Dindent{\th@Dheadingshape ##1.}\Dquad     
        \D@restoremargins}%
  \def\@opargbegintheorem##1##2##3{\def\next{##3}%
  \item[]\th@Dindent{\th@Dheadingshape ##1\ifx\next\empty
  \else\ ##3\fi.}\Dquad         
        \D@restoremargins}}

\gdef\th@DthAndSuchliketheo{
  \D@makemargins%
  \def\@begintheorem##1##2{%
    \@latex@error{likethm: You must provide an argument in square brackets,
    though it may be empty [] !}%
    }%
  \def\@opargbegintheorem##1##2##3{%
        \def\next{##3}\ifx\next\empty\item[\th@Dindent]\else
        \item[]\th@Dindent{\th@Dheadingshape \next.}\Dquad\fi
        \D@restoremargins}}

\theoremstyle{DthAndSuchliketheo}
\theorembodyfont{\bigpenalty\itshape}   

\theorembodyfont{\rmfamily}  

\theoremstyle{DthAndSuchtheostar}
\theorembodyfont{\bigpenalty\itshape}
\newtheorem{thm*}{\theoname}
\newtheorem{lem*}{\lemmaname}
\newtheorem{pro*}{\propositionname}
\newtheorem{cor*}{\corollaryname}
\newtheorem{conjecture*}{\conjecturename}
\theorembodyfont{\smallpenalty\rmfamily}
\newtheorem{notation*}{\notationname}
\newtheorem{exa*}{\examplename}
\newtheorem{examples*}{\examplesname}
\newtheorem{definition*}{\definitionname}
\newtheorem{definitions*}{\definitionsname}
\newtheorem{rem*}{\remarkname}
\newtheorem{remarks*}{\remarksname}
\theoremstyle{DthAndSuchtheo}
\theorembodyfont{\bigpenalty\itshape}
\newtheorem{thm}{\theoname}[section]
\newtheorem{lem}[thm]{\lemmaname}
\newtheorem{pro}[thm]{\propositionname}
\newtheorem{cor}[thm]{\corollaryname}

\theorembodyfont{\smallpenalty\rmfamily}
\newtheorem{exa}[thm]{\examplename}

\newtheorem{rem}[thm]{\remarkname}
\newtheorem{remarks}[thm]{\remarksname}

\theoremstyle{DthAndSuchtheo}               
\theorembodyfont{\itshape}

\newcommand{\proof}[1][]{{\par\removelastskip\vskip.6\baselineskip   
    \noindent\th@Dindent\def\next{#1}%
    {\itshape\proofname\ifx\next\empty\else\next\fi\ifperiod.%
      \else\global\periodtrue\fi\Dquad}%
    \clubpenalty=5000\rm\ignorespaces}\setcounter{step}{0}\setcounter{case}{0}}

\newcounter{step}
\newcommand{\step}{\stepcounter{step}
  \par\indent{\itshape Step \thestep .\hspace{1ex}}}

\newcounter{case}

\def\case{\stepcounter{case}%
  \vspace{.5\baselineskip}
  \def\mynext{\thecase}
  \ifnum\mynext=1{\itshape First case.}
  \else
    \ifnum\mynext=2{\itshape Second case.}
    \else
      \ifnum\mynext=3{\itshape Third case.}
      \else
        \ifnum\mynext=4{\itshape Fourth case.}
        \fi      
      \fi
    \fi
  \fi
  \hspace{1ex}}

\newcommand{\likeproof}[1][]{{\par\removelastskip\vskip.6\baselineskip
    \noindent\th@Dindent\def\next{#1}%
    {\itshape\ifx\next\empty\else\next\fi\ifperiod.%
      \else\global\periodtrue\fi\Dquad}%
    \clubpenalty=5000\rm\ignorespaces}\hspace{-2pt}\setcounter{step}{0}}

\def\qed{{\ifmmode\hskip 6mm plus 1mm minus 3mm{$\square$}
    \else
    \hfil\hskip1em\null\nobreak\hfil
    {\clubpenalty=50000\bigpenalty\hfill
      $\square$\parfillskip=0pt\finalhyphendemerits=0
      \let\par=\endgraf\par}
    \fi
    \Dpar\vskip.6\normalbaselineskip}}

\makeatother

\newcommand{\typeI}{type I}
\newcommand{\typeII}{type II}
\newcommand{\reduced}[1]{{#1}^{\text{\normalfont red}}}
\newcommand{\Irr}{\operatorname{Irr}}

\numberwithin{equation}{section}

\begin{document}

\title{Twisted Kodaira-Spencer classes and the geometry of surfaces of
general type}
\author{Daniel Naie and Igor Reider}


\begin{flushright}
  {\it To Fedor Bogomolov on his 65th birthday}  
\end{flushright}

\bigskip

\begin{abstract}
We study the cohomology groups $H^1(X,\Theta_X(-mK_X))$, for $m\geq1$,
where $X$ is a smooth minimal complex surface of general type,
$\Theta_X$ its holomorphic tangent bundle, and $K_X$ its canonical
divisor.  One of the main results is a precise vanishing criterion for
$H^1(X,\Theta_X (-K_X))$ (Theorem \ref{t:theResult}).

The proof is based on the geometric interpretation of non-zero
cohomology classes of $H^1(X,\Theta_X (-K_X))$. This interpretation in
turn uses higher rank vector bundles on $X$.

We apply our methods to the long standing conjecture saying that the
irregularity of surfaces in $\PP^4$ is at most $2$.  We show that if
$X$ has prescribed Chern numbers, no irrational
pencil, and is embedded in $\PP^4$ with a sufficiently large degree,
then the irregularity of $X$ is at most $3$.

\end{abstract}

\bigskip


\section{Introduction}

From Kodaira-Spencer theory of deformations of complex structures,
\cite{KS}, it is known that for a smooth complex projective variety
$X$ with holomorphic tangent bundle $\Theta_X$, the cohomology group
$H^1(X,\Theta_X)$ parametrizes the infinitesimal deformations of
complex structures on $X$---hence the importance to study
this group. In this paper we propose to study its twisted versions.
More precisely, let $X$ be a smooth minimal complex surface of general
type and consider the cohomology groups $H^1(X,\Theta_X(-m K_X))$,
where $m$ is a positive integer and $K_X$ the canonical divisor of
$X$.  The basic question we wish to address is the vanishing
(respectively the non-vanishing) of these groups.  It is customary to
call elements of $H^1(X,\Theta_X)$ Kodaira-Spencer classes. Following
this tradition, we will refer to classes of $H^1(X,\Theta_X(-m K_X))$
as twisted Kodaira-Spencer classes of degree $m$.

Of course, if $K_X$ is ample the general result of Serre on sheaf
cohomology implies the vanishing of $H^1(X,\Theta_X(-m K_X))$ for all
$m$ sufficiently big.  On the other hand, if $K_X$ is not ample, then
$X$ carries $(-2)$ curves, \ie smooth rational curves $C$ with 
$C^2 = -2$ and $C\cdot K_X =0$.  Furthermore, it is easy to see that
every $(-2)$ curve on $X$ gives rise to a non-zero class in
$H^1(X,\Theta_X(-m K_X))$, for {\it all} positive $m$ (see (b) in 
Remark~\ref{r:extensions}).  This observation
indicates that the twisted Kodaira-Spencer classes contain a rather
precise geometric data. In fact our main line of investigation as well
as the strategy for proving the vanishing are very much modeled on
this observation. Namely, we seek to extract geometric data from
non-zero classes in $H^1(X,\Theta_X(-m K_X))$. Once this is done, we
impose suitable hypotheses on $X$ to rule out the existence of such
data and thus obtain the vanishing of $H^1(X,\Theta_X(-m K_X))$. The
following is one of the main results of the paper illustrating our
approach.

\begin{thm} \label{t:theResult}
Let $X$ be a smooth minimal complex surface of general type such that
its Chern numbers satisfy $c_2(X)/K^2_X<{5}/{6}$.  If $X$ contains no
smooth rational curves with self-intersection $-2$ or $-3$, then
\[
  H^1(X,\Theta_X (-K_X))=0.
\] 
\end{thm}

This result points to the following general paradigm regarding the
groups\linebreak $H^1(X,\Theta_X(-mK_X))$, for $m\geq1$:

\begin{center}
\begin{minipage}[c]{.92\linewidth}
{\it If the topological invariants of $X$ are suitably constrained,
then the \mbox{non-vanishing} of $H^1(X,\Theta_X(-mK_X))$ implies the
existence of configurations of special curves on $X$.} 
\end{minipage}  
\end{center}

This heuristic principle is realized in case $m=1$ by the following
technical result which constitutes the backbone of our considerations
and might be of independent interest.

\begin{pro}
  \label{p:decomp0}
Let $X$ be a smooth minimal complex surface of general type such that
its Chern numbers satisfy $c_2(X)/K^2_X <1$. If
$H^1(X,\Theta_X(-K_X))\neq 0$, then the canonical divisor $K_X$ admits
a decomposition
\begin{equation}
  \label{f:decomp}
  K_X = L + E
\end{equation}
subject to the following conditions:
\begin{enumerate}\listspace
\item[(i)]
  $L$ has the Iitaka dimension $2$, and
\item[(ii)]
  $E$ is a non-zero effective divisor which is numerically effective
(nef), unless $E$ contains irreducible components $C$ such that
$C\cdot E =-2$ or $-1$.  Such components are smooth rational curves of
self-intersection $-2$ or $-3$.
\end{enumerate}
\end{pro}

A way to prove Theorem \ref{t:theResult} can now be summarized as
follows: The hypotheses of Theorem \ref{t:theResult} impose severe
restrictions on the decomposition in \eqref{f:decomp} and these
restrictions, in their turn, lead to a contradiction.  Hence
$H^1(X,\Theta_X(-K_X))$ must vanish.  

\bigskip

Now, we will explain how to pass from the non-vanishing of
$H^1(X,\Theta_X(-K_X))$ to the decomposition in \eqref{f:decomp} and
to its various geometric properties.  The main idea is to use the
identification
\[
H^1(X,\Theta_X(-K_X)) \cong \Ext^1(\Omega_X,\Oo_X(-K_X)),
\]
where $\Omega_X$ is the cotangent bundle of $X$ (see \cite{R} for a
similar approach).  We can and will view a non-zero class 
$\xi\in H^1(X,\Theta_X(-K_X))$ as the corresponding extension, \ie the
short exact sequence of locally free sheaves on $X$
\begin{equation}
  \label{eq:ext-seq0}
  0 \lra \Oo_X(-K_X) \lra \Tt_\xi \lra \Omega_X \lra 0.
\end{equation}
Our attention shifts to the middle term sheaf $\Tt_\xi$.  It is locally
free of rank $3$, with Chern invariants
\[
c_1(\Tt_\xi)=0 \textq{and} c_2(\Tt_\xi)=c_2-K_X^2.
\]
The hypothesis $c_2(X)/K^2_X <1$ yields that $\Tt_\xi$ is unstable in
the sense of Bogomolov.  Considering its {\it Bogomolov destabilizing
subsheaf}\/ is the key to obtaining the decomposition in
\eqref{f:decomp}.  This subsheaf, which we call $\Ff$, must be of rank
$2$, and its determinant provides the divisor $L$ in the decomposition
\eqref{f:decomp}.  To see the other part of the decomposition, we
observe that the inclusion of $\Ff$ in $\Tt_\xi$ combined with the
epimorphism in \eqref{eq:ext-seq0} induces a morphism
\[
\Ff \longrightarrow \Omega_X
\]
which is generically an isomorphism. Furthermore its rank must drop
along a non-zero divisor.  This is the divisor $E$ in
\eqref{f:decomp}.  With the decomposition established, it is rather
straightforward to obtain some information on the numerical properties
of the decomposition, such as $L^2$, $L\cdot K_X$, $E^2$, $E\cdot K_X$
(see \eqref{eq:mainInequality}.  But one can go further on by getting
rather detailed information concerning the irreducible components of
$E$.  This constitutes the most technical part of the paper.

The above discussion clearly shows that the non-zero twisted
Kodaira-Spencer classes contain quite precise geometric information
and one might be tempted to say that for understanding the
geometry of $X$, it is more useful to have the non-vanishing of the
groups $H^1(X,\Theta_X(-mK_X))$ for some $m>0$, rather than
their vanishing.

Theorem \ref{t:theResult} belongs to the category of ``effective
vanishing theorems".  It could also be viewed as a cohomological
criterion for the existence of particular rational curves on surfaces
of general type.  This, in itself, can be quite helpful.  At the same
time, the groups of the more general form $H^1(X,\Theta_X(-D))$, with
$D$ a ``positive" divisor, come up in various contexts.  So the
control of these groups can be useful for applications.

One immediate application comes from the deformation theoretic
interpretation of $H^1(X,\Theta_X(-D))$.  Similar to the usual
Kodaira-Spencer classes, the twisted ones in $H^1(X,\Theta_X(-D))$ can
also be viewed as infinitesimal deformations.  Namely, one considers
deformations of pairs $(X,C)$, where $X$ is a surface of general type
and $C$ is a smooth curve in the linear system $\mid D\mid$.  The
group $H^1(X,\Theta_X(-D))$ parametrizes the infinitesimal
deformations of $(X,C)$ where the corresponding infinitesimal
deformation of the curve $C$ is trivial. In other words
$H^1(X,\Theta_X (-D))$ can be identified with the Zariski tangent
space along the fibre through $(X,C)$ of the forgetful functor
\begin{equation*}
  F : (X,C) \mapsto C.
\end{equation*}
So the vanishing of $H^1(X,\Theta_X(-D))$ would tell us that this
functor gives an immersion between the corresponding stacks. Thus
Theorem \ref{t:theResult} yields the immersion of $\MMM^p_{K^2,c_2}$,
the stack of pairs $(X,C)$, where $C$ is a smooth curve in the
canonical linear system $|K_X|$ and $X$ is subject to the hypotheses
of Theorem \ref{t:theResult}, into $\MMM_{K^2+1}$, the stack of smooth
projective curves of genus $g_C=K^2_X+1$.

\begin{cor}
  \label{c:immer}
If\/ $\dfrac{c_2}{K^2}<\dfrac{5}{6}$, then the forgetful functor
$F : \MMM^p_{K^2,c_2} \to \MMM_{K^2+1}$
is an immersion if the canonical system $|K_X|$ contains a smooth
curve.
\end{cor}

\bigskip

The second application developed in the paper comes from the general
observation that the groups of the form $H^1(X,\Theta_X(-D))$ ($D$ is,
as before, some positive divisor on $X$) control a certain amount of
the {\it extrinsic} geometry of $X$.  More precisely, if $X$ is
embedded in a smooth projective variety $Y$, then the normal bundle
$\Nn_{X/Y}$ of $X$ in $Y$ is related to the tangent bundle $\Theta_X$
of $X$ via the normal sequence.  Hence the relation between
$H^0(X,\Nn_{X/Y}(-D))$ and $H^1(X,\Theta_X(-D))$.

It may happen that for some {\it a priori} ``extrinsic" reasons, one
knows that the cohomology group $H^0(X,\Nn_{X/Y}(-D))\neq 0$.  Then,
provided that the coboundary map
\[
H^0(X,\Nn_{X/Y}(-D)) \lra H^1(X,\Theta_X(-D))
\]
is not zero, one obtains the non-vanishing of $H^1(X,\Theta_X(-D))$.
This, according to our heuristic principle, should impose topological
and geometrical constraints on $X$.  It potentially opens a way to
having restrictions on the topology/geometry of surfaces embeddable in a
given smooth projective variety $Y$.  We summarize this ``extrinsic --
intrinsic" reasoning in the following diagram.

\bigskip

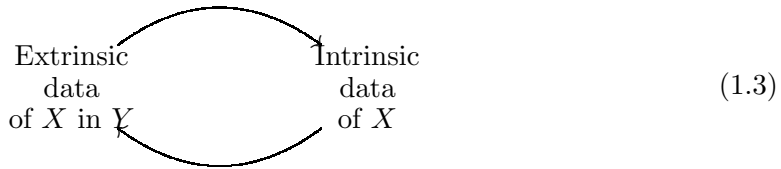
\begin{equation}
  \label{d:e-i}
  \UseTips
  \xymatrix @R=200pt @W=5mm @H=7mm @*[r]@*[u]{
    *\txt{Extrinsic\\ data\\ of $X$ in $Y$}
    \ar@/^2.5pc/[rr] &&
    *\txt{Intrinsic\\ data\\ of $X$} \ar@/^2.5pc/[ll] }
\end{equation}

\bigskip

We apply this line of thinking to surfaces in $\PP^4$ and in
particular to the long standing conjecture about the upper bound on
the irregularity of surfaces in $\PP^4$. This says that the
irregularity of such surfaces is at most $2$.  For a more ample
discussion of the history and the background of this problem we refer
the reader to \cite{OSS} and the references therein.  To our
understanding, there is no conceptual reason for such an upper bound
and this estimate is largely based on the lack of examples. However,
by exploiting non-zero twisted Kodaira-Spencer classes for appropriate
twists of $\Theta_X$, we are able to show the following.

\begin{thm}[= Theorem~\ref{bound-q}]
  \label{t0:bound-q}
Given a positive integer $\chi$ and a rational $\beta$,
$2\leq \beta\leq9$, there is a number $d_0(\chi,\beta)$ such every smooth
minimal surface $X\subset\PP^4$ of degree $d>d_0(\chi,\beta)$, with
$\chi(\Oo_X)=\chi$, $K_X^2=\beta\chi(\Oo_X)$, and having no irrational
pencil, has the irregularity at most $3$.  Furthermore, if the
irregularity is equal to $3$, then  $X$ must be subject to the condition 
(i) or (ii) in Lemma \ref{l:4-5}.
\end{thm}

Let us explain how this result and its proof fit into the general
heuristic scheme \eqref{d:e-i}.  To begin with, the constant
$d_0(\chi,\beta)$ and the condition $d>d_0(\chi,\beta)$ of
Theorem~\ref{t0:bound-q} provide an {\it a priori} extrinsic piece of
information for the pair $(X,\PP^4)$ on the left side of the diagram
\eqref{d:e-i}.  This datum is the existence of a $3$-fold of degree
$m_X\leq 5$ in $\PP^4$ containing $X$.  This observation is
essentially due to Ellingsrud and Peskine in \cite{E-P}.  Then, the
right oriented arrow in \eqref{d:e-i} translates the existence of this
$3$-fold into the intrinsic datum of the non-vanishing of
$H^1(X,\Theta_X(-D))$, where $D=K_X+(5-m_X)H$, and where $H$ denotes
the divisor class of hyperplane sections of $X \subset \PP^4$.  A
version of Proposition \ref{p:decomp0} implies a decomposition, as in
\eqref{f:decomp}, of the canonical divisor.  It should be pointed out
that in this version we use the hypothesis on irrational pencils
instead of the inequality between the Chern numbers of Proposition
\ref{p:decomp0}.  The intrinsic datum of the decomposition of $K_X$ is
then translated back into extrinsic data by establishing a
relationship between $L$ in \eqref{f:decomp} and the hyperplane
section $H$---this is the left oriented arrow in \eqref{d:e-i}.  Our
bound on the irregularity of $X$ is an immediate consequence of this
relation.

Though Theorem \ref{t0:bound-q} is a long way from the unconditional
bound of the conjecture, it seems to be useful since it
suggests that one of the conceptual reasons for not having surfaces in
$\PP^4$ with large irregularity is the {\it non-vanishing} of
cohomology groups of the form $H^1(X,\Theta_X(-D))$ for an
appropriately chosen ``positive" divisor $D$ on $X$.  

The arguments developed in the proof of Theorem~\ref{t0:bound-q} can
also be used to show that for high degree surfaces in $\PP^4$ with
bounded holomorphic Euler characteristic, their topological index is
negative.  (See Theorem~\ref{t:topCondition} for the precise
statement.)

\bigskip

Let us end this section with some comments and questions.
Theorem \ref{t:theResult} together with the additional assumption that
the canonical linear system $\mid K_X \mid$ contains a reduced
irreducible member implies the vanishing of $H^1(X,\Theta_X(-mK_X))$
for all $m\geq1$ (Corollary \ref{c:vanish}). This gives a complete
answer to the question about the vanishing of the groups
$H^1(X,\Theta_X(-mK_X))$, $m\geq1$, for a large class of surfaces of
general type.  However, the situation is not fully satisfactory with
regard to the hypothesis in Theorem \ref{t:theResult} concerning the
$(-3)$-curves---we do not know if this assumption is really necessary.
Let us explain this point. Contrary to $(-2)$-curves, the
$(-3)$-curves, by themselves, do not produce cohomology classes in
$H^1(X,\Theta_X(-K_X))$. These curves appear naturally in the course
of the study of the decomposition of $K_X =L+E$ in
\eqref{f:decomp}. But each $(-3)$-curve must appear in the divisor $E$
as a part of a rather involved configuration of other curves, and at
this stage, we are unable to treat these configurations efficiently.

It should also be noticed that the problem of $(-3)$-curves
disappears for the groups $H^1(X,\Theta_X(-mK_X))$, with $m\geq 4 $, and
our approach allows to have a purely cohomological criterion for the
existence of $(-2)$-curves on surfaces of general type (that will
appear elsewhere).

We also believe that the groups $H^1(X,\Theta_X(-mK_X))$, for all $m$
sufficiently large, must be spanned by the twisted Kodaira-Spencer
classes associated to $(-2)$-curves. Another way to put it, the
twisted Kodaira-Spencer classes of sufficiently high degree should
provide a stronger version of the decomposition of type
\eqref{f:decomp} where the divisor $E$ is composed entirely of
$(-2)$-curves.

We limited our discussion to surfaces of general type.  However, it is
clear that our approach is valid in higher dimensions as well and one
can ask if it is possible to have a result analogous to
Theorem~\ref{t:theResult} in dimensions bigger than $2$.

\bigskip

The paper is organized as follows.  In \S2, we prove the decomposition
\eqref{f:decomp} and derive the basic properties of the divisor $L$ in
this decomposition (see Lemma \ref{l:decomp} and Lemma \ref{l:L2-LK}).
The third section is devoted to a detailed study of the irreducible
components of the divisor $E$ in \eqref{f:decomp}.  In \S4, we give a
proof of Theorem \ref{t:theResult}.  Finally, in \S5, the two results
about surfaces in $\PP^4$ are considered (see Theorem
\ref{bound-q} and Theorem~\ref{t:topCondition}).

\bigskip

\paragraph{Notation and conventions}

$X$ is a smooth minimal complex surface of general type unless
otherwise stated.

$K_X$ is the canonical divisor of $X$.

$K^2_X$, $c_2(X)$, and $\chi(\Oo_X)$ are the Chern numbers and the
holomorphic Euler characteristic of $X$ respectively.

$\alpha_X:=\dfrac{c_2(X)}{K_X^2}$.

$\Theta_X$ is the holomorphic tangent bundle of $X$ and $\Omega_X$ is
its holomorphic cotangent bundle.
$p_g(X)=h^2(X,\Oo_X)=h^0(X,\Oo_X(K_X))$ and
$q(X)=h^1(X,\Oo_X)=h^0(X,\Omega_X)$ are respectively, the geometric
genus and the irregularity of $X$.

For a coherent sheaf $\Gg$ on a variety $M$, the cohomology group
$H^i(M,\Gg)$ is denoted by $H^i(\Gg)$ if no confusion is likely. 
  
For an effective divisor $D=\sum_\alpha \mu_\alpha C_\alpha$ on a
surface $X$, where the $C_\alpha$ are reduced irreducible curves,
$\reduced{D}$ denotes the reduced associated divisor 
$\sum_\alpha C_\alpha$.

All equalities between divisors take place in $A^1(X)$ of the Chow
ring of $X$ except the ones involving the Zariski decomposition.  In
these cases, the equalities are considered in the N\'eron-Severi group
$NS(X)_\QQ=NS(X)\otimes_\ZZ\QQ$.


\section{The non-vanishing of $H^1(X,\Theta_X(-K_X))$ and \\
  the decomposition of $K_X$}

In this section we show how the existence of a non-trivial element 
in $H^1(X,\Theta_X(-K_X))$ gives rise to a distinguished
decomposition of the canonical divisor $K_X$ of $X$.

\begin{lem}
  \label{l:decomp}
Assume the Chern numbers $K^2$ and $c_2(X)$ are subject to
\[
  \alpha_X=\dfrac{c_2(X)}{K_X^2}<1 
\]
and assume $H^1 (X,\Theta_X (-K_X))$ to be non-zero.  Then the
canonical divisor $K_X$ admits the decomposition
\[
  K_X = L+E
\]
where $L$ has Iitaka dimension $2$ and $E$ is effective and non-zero. 
\end{lem}

\proof
A non-zero element $\xi \in H^1(X,\Theta_X (-K_X))$ via the natural
identification
\[
   H^1(X,\Theta_X (-K_X)) \cong \Ext^1(\Omega_X,\Oo_X(-K_X))
\]
defines the extension 
\begin{equation*}
  \label{eq:ext}
  0 \lra \Oo_X(-K_X) \lra \Tt_\xi \lra \Omega_X \lra 0.
\end{equation*}
The sheaf $\Tt_\xi$, seating in the middle of the sequence, is locally
free of rank $3$ and its Chern invariants are
\[
c_1(\Tt_\xi)=0 \textq{and} c_2(\Tt_\xi)=c_2-K_X^2.
\] 
The assumption $\alpha_X<1$ implies that $\Tt_\xi$ is Bogomolov
unstable.  Let $\Ff$ be its Bogomolov destabilizing subsheaf.
One knows that its determinant $\det\Ff=\Oo_X(L)$ is in the positive
cone of $NS(X)_\QQ$ (see \cite{Bog}).  This and the fact
that the cotangent sheaf $\Omega_X$ cannot have subsheaves of rank
$1$ and of Iitaka dimension $2$ imply that the rank of $\Ff$ is $2$.
Furthermore, we may assume it to be saturated.  Then the quotient
$\Tt_\xi/\Ff$ is torsion-free and $\Tt_\xi/\Ff\simeq\Jj_Z(-L)$, where
$Z$ is a subscheme of codimension $2$ and $\Jj_Z$ is its sheaf of
ideals.  This yields the following diagram.
\begin{equation}
  \label{d:unstable}
  \UseTips
  \xymatrix@W=5mm@H=5mm{
    && 0 \ar[d] \\
    && \Ff \ar[d]\ar[rd]^{\phi_\xi} \\
    0 \ar[r] & \Oo_X(-K_X) \ar[r] \ar[rd]_{\psi_\xi}&
    \Tt_\xi \ar[r]\ar[d] & \Omega_X \ar[r] & 0 \\
    && \Jj_Z(-L) \ar[d] \\
    && 0
  }
\end{equation}
Using again the fact that $\Omega_X$ has no rank $1$ subsheaves of
Iitaka dimension $2$, we conclude that the morphism $\phi_\xi$ in
\eqref{d:unstable} is generically an isomorphism.  In addition, its
cokernel, $\coker{\phi_\xi}$, is a sheaf supported on a codimension
$1$ subscheme since otherwise $\phi_\xi$ would be an isomorphism
making the extension class $\xi$ trivial.  The resulting exact
sequence
\[
0 \lra \Ff \stackrel{\phi_\xi}{\lra} \Omega_X \lra 
\coker{\phi_\xi} \lra 0 
\]
yields the asserted decomposition
\[
  K_X=L+E
\]
with $L$ of Iitaka dimension $2$ and $E=c_1 (\coker{\phi_\xi})$
effective and non-zero.
\qed

We begin the study of the decomposition of $K_X$ in Lemma
\ref{l:decomp} by recording some properties of $L$.

\begin{lem}
  \label{l:L2-LK}
The divisor $L$ in Lemma \ref{l:decomp} admits the Zariski
decomposition%
\footnote{As we have specified in the introduction, the equality
  $L=L^++L^-$ takes place in $NS(X)_\QQ$.}
\[
  L = L^+ + L^-
\]
where $L^+$, the positive part of $L$, is nef and big, and $L^-$ is
the negative part of $L$. Furthermore, $L$ is subject to the following
numerical conditions:
\begin{equation}
  \label{eq:mainInequality}
  L^2 \geq \frac{3}{2}\,(1-\alpha_X)\,K_X^2
  \textq{and}
  L\cdot K_X \geq \sqrt{\frac{3}{2}(1-\alpha_X)}\,K_X^2.
\end{equation}
\end{lem}
 
\proof
By construction, $L$ has Iitaka dimension $2$, hence it has the
Zariski decomposition $L=L^++L^-$ as asserted.

The two inequalities of the lemma come from the diagram
\eqref{d:unstable}.  Namely, combining its vertical and horizontal
sequences yields
\[
  c_2(X)-K_X^2 = c_2(\Ff)-L^2+\deg Z.
\]
Since $\Ff$ is a rank $2$ subsheaf of the cotangent sheaf $\Omega_X$
and $\det(\Ff)=\Oo_X(L)$ has Zariski decomposition, a
result of Miyaoka, \cite[Remark~4.18]{Miy}, says that its Chern
invariants are subject to
\[
  c_2(\Ff) \geq \frac{1}{3}\,L^2-\frac{1}{12}\,(L^-)^2.
\]
Combining the two inequalities, we obtain
\begin{equation}
  \label{eq:LSquareAndKSquare}
  L^2 \geq 
  \frac{3}{2}(K_X^2-c_2(X))+\frac{3}{2}\deg Z - \frac{1}{8}(L^-)^2
  \geq \frac{3}{2}(1-\alpha_X)K_X^2,
\end{equation}
thus the first of the asserted inequalities. The second follows from
the first and the Hodge index theorem 
$(L\cdot K_X)^2\geq K^2L^2$.
\qed

The following lemma compares the two divisors, $L$ and $E$, of the
decomposition in Lemma~\ref{l:decomp}.

\begin{lem}
  \label{l:positiveCone}
Let $K_X=L+E$ be the decomposition in Lemma~\ref{l:decomp}.

{\sc 1)} If $\alpha_X<5/6$, then $L\cdot K_X>E\cdot K_X$.

{\sc 2)} If $\alpha_X\leq1/2$, then $(L-E)^2\geq0$. 
\end{lem}

\proof
To prove 1), assume\/ $L\cdot K_X\leq E\cdot K_X$.  Then,
\[
L\cdot K_X \leq E\cdot K_X \leq K_X^2-L\cdot K_X,
\]
and using
\eqref{eq:mainInequality}, we obtain
\[
K_X^2 \geq 2L\cdot K_X \geq 2\sqrt{\frac{3}{2}(1-\alpha_X)}\,K_X^2.
\]
Hence 
$1\geq \sqrt{6(1-\alpha_X)}$,
\ie $\alpha_X\geq5/6$.

\medskip

To prove 2) we use the hypothesis $\alpha_X\leq1/2$ and
\eqref{eq:mainInequality} to arrive at
\begin{multline*}
  (L-E)^2 = (2L-K_X)^2 
  = 4L^2-4L\cdot K_X+K_X^2
  \geq 6(1-\alpha_X)K_X^2-4L\cdot K_X+K_X^2 \\
  \geq 3K_X^2-4L\cdot K_X+K_X^2
  = 4(K_X^2-L\cdot K_X) = 4 E\cdot K_X \geq 0.
\end{multline*}
\qed

In the sequel, we will need the following cohomological property 
of the divisor $E$.

\begin{lem}
  \label{l:E-coh}
Let $e$ be a section of $\Oo_X (E)$ corresponding to the divisor $E$.
Then $e\cdot\xi=0$ in $H^1(\Theta_X(E-K_X))$.
\end{lem}

\proof
Dualizing the diagram \eqref{d:unstable} and tensoring it with
$\Oo_X(-L)$ give the following diagram.
\begin{equation}
  \label{d:unstableD}
  \UseTips
  \xymatrix@W=5mm@H=5mm{
    && \Oo_X \ar[d]\ar[rd]^e \\
    0 \ar[r] & \Theta_X(-L) \ar[r] &
    \Tt_\xi^\ast(-L) \ar[r] & \Oo_X(K_X-L) \ar[r] & 0 \\
  }
\end{equation}
We see that the coboundary map 
\[
  H^0(\Oo_X(E)) \lra H^1(\Theta_X(-L)) = H^1(\Theta_X(E-K_X))
\]
given by the cup-product with the class $\xi$ contains the section $e$
in its kernel.
\qed

\begin{rem}
  \label{r:extensions}
(a) The equality $e\cdot\xi=0$ in Lemma~\ref{l:E-coh} implies that $\xi$
is supported on the divisor $E=\{e=0\}$.  More precisely, from the
exact sequence
\[
  0 \lra H^0(\Theta_X\otimes\Oo_E(E-K_X)) \lra 
  H^1(\Theta_X(-K_X)) \stackrel{e}{\lra}
  H^1(\Theta_X(E-K_X))
\]
it follows that $\xi$ is the image of a unique global section of
$\Theta_X\otimes\Oo_E(E-K_X)$.

(b) In case $E$ is a $(-2)$-curve, one has a converse, \ie $E$ defines
a unique (up to a non-zero scalar multiple) non-zero cohomology class
in $H^1(\Theta_X(-K_X))$.  This follows since
\[
  \Theta_X\otimes\Oo_E(E-K_X) = \Theta_X\otimes\Oo_E(E)
  \cong \Theta_X\otimes\Oo_{\PP^1}(-2) \cong 
  \Oo_{\PP^1}\oplus\Oo_{\PP^1}(-4),
\]
and hence 
$H^0(\Theta_X\otimes\Oo_E(E-K_X))=H^0(\Oo_{\PP^1}\oplus\Oo_{\PP^1}(-4))=\CC$. 
The same argument is valid for all $H^1(\Theta_X(-mK_X))$, where
$m\geq1$.  Thus, a $(-2)$-curve gives rise to a non-zero cohomology
class in every group $H^1(\Theta_X(-mK_X))$, $m\geq1$.
\end{rem}

\section{The study of the irreducible components of $E$}

This section is entirely devoted to a study of the divisor $E$ of the
decomposition in Lemma~\ref{l:decomp}.

\begin{lem}
The divisor $E$ passes through the subscheme $Z$.
\end{lem}

\proof
To see this, it suffices to show that $h^0(X,\Jj_Z(E))\neq 0$.  This
follows by tensoring the morphism $\psi_\xi$ in \eqref{d:unstable}
with $\Oo_X(K_X)$.
\qed

From the vertical sequence in \eqref{d:unstable}, the subscheme $Z$ can
be seen as the zero-locus of a section 
$s_\xi\in H^0(X,\Tt^\ast_\xi(-L))$.  Using this section, we can define
the restriction of $Z$ to a component of $E$.  More precisely, for
every component $C$ of $E$, denote by $Z_C=\{s_\xi|_C=0\}$.

The main tool for a detailed study of $E$ is the following construction.
We take the second exterior product in diagram \eqref{d:unstable} to
obtain the diagram,
\begin{equation}
  \label{d:exteriorProduct}
  \UseTips
  \xymatrix@W=5mm@H=5mm{
    && \Oo_X(L) \ar[d]\ar[rd]^e \\
    0 \ar[r] & \Theta_X \ar[r] &
    \bigwedge^2\Tt_\xi\ar[r] & \Oo_X(K_X) \ar[r] & 0 
  }
\end{equation}
where $e$ is a section of $\Oo_X (E)$ corresponding to the divisor
$E$.  For any component $C$ of $E$, the morphism in
\eqref{d:exteriorProduct} given by $e$ vanishes over $C$ and hence it
induces a non-zero morphism
\[
  \tau_C:\Oo_C(L) \to \Theta_X\otimes\Oo_C
\]
whose zero locus is $Z_C$.  If $C$ is a reduced and irreducible
component of $E$, we combine this morphism with the normal sequence of
$C$ in $X$ to obtain the following commutative diagram.
\[
\UseTips
\xymatrix@W=5mm@H=5mm{
  && \Oo_C(L) \ar[d]_{\tau_C} \ar[rd]^{\psi_C} \\
  0 \ar[r] & \Theta_C \ar[r] &
  \Theta_X\otimes\Oo_C \ar[r] & \Oo_C(C) 
}
\]
The irreducible components $C$ of $E$ will be distinguished
according to the properties of the morphism $\psi_C$.

\smallskip

{\bf\typeI}: If $\psi_C=0$, in this case $\tau_C$ factors through
the tangent sheaf $\Theta_C$ of $C$ and the degree of $Z_C$ is subject
to
\begin{equation}
  \label{eq:typeI}
  \deg Z_C \leq \deg \eta_C^\ast(Z_C) = 2-2g(\Ctilde)-C\cdot L,
\end{equation}
where $\eta_C:\Ctilde\to C$ is the normalization of $C$ (see \cite{R},
pp.432-433, for details). Such a component will be called of \typeI.

\smallskip

{\bf\typeII}:  If $\psi_C\neq0$, in this case $\psi_C$ defines a
non-zero section of $\Oo_C(C-L)$ vanishing on $Z_C$, and hence
\begin{equation}
  \label{eq:typeII}
  \deg Z_C \leq C^2-C\cdot L.
\end{equation}
We will call such a component of \typeII.

\bigskip

The divisor $E$ can be written as the sum of two parts,
\begin{equation}
  \label{ty-d}
  E = E_I+E_{II},
\end{equation}
with
\[
  E_I := \sum_{C\text{ of \typeI}}m_C\,C
  \textq{and}
  E_{II} := \sum_{C\text{ of \typeII}}m_C\,C,
\]
where $m_C$ denotes the multiplicity of a component $C$ in $E$.  

The following lemma relates the components of $E_{II}$ and the Zariski
decomposition of $L$ in Lemma \ref{l:L2-LK} under the additional
assumption $\alpha_X < \dfrac{5}{6}$ of Theorem \ref{t:theResult}.

\begin{lem}
  \label{l:typeIIandLminus}
Let $\alpha_X < \dfrac{5}{6}$. Then every irreducible component $C$ of
$E$ of type II must also be a component of $L^-$ and
\[
  C\cdot L \leq C^2 < 0.
\]
Furthermore, the multiplicity of $C$ in $L^-$ is $\geq1$.
\end{lem}

\proof
From \eqref{eq:typeII}, we deduce $C\cdot L\leq C^2$.  First, we claim 
that
$C\cdot L\leq0$.  Indeed, if $C\cdot L$ is positive, then $C^2>0$.
From this and the Hodge index theorem, $(C\cdot L)^2\geq C^2\,L^2$, it
follows that
\begin{equation}
  \label{eq:anIneq}
  C^2 \geq C\cdot L \geq L^2   \geq \frac{3}{2}\,(1-\alpha_X) K_X^2
\end{equation}
where the last inequality comes from Lemma \ref{l:L2-LK}.  This and
the Hodge index $(C\cdot K_X)^2 \geq C^2\,K_X^2$ yield
\[
E\cdot K_X \geq C\cdot K_X
\geq \sqrt{\frac{3}{2}\,(1-\alpha_X)} K_X^2.
\]
Combining this with the second inequality of Lemma \ref{l:L2-LK}
we obtain
\[
K_X^2 = L\cdot K_X + E\cdot K_X
\geq 2\,\sqrt{\frac{3}{2}\,(1-\alpha_X)} K_X^2. 
\]
But this yields $\alpha_X\geq5/6$, contrary to our assumption.

Next combining the inequality $C\cdot L\leq 0$ with the Zariski
decomposition $L=L^+ + L^-$ we obtain that $C^2 < 0$.  This, together
with $L\cdot C\leq C^2$, yields $L\cdot C<0$.  Hence $C$ must be a
component of $L^-$.

Turning to the last assertion of the lemma, let $\mu_C$ be the
multiplicity of $C$ in $L^-$.  Then
\[
  0 \geq C\cdot L-C^2 
  = C\cdot (L^++L^-)-C^2 
  =C\cdot L^--C^2
  \geq (\mu_C-1) C^2,
\] 
hence $\mu_C\geq 1$.
\qed

In analyzing the components of \typeI, we will repeatedly make use of
the following lemma. 
 
\begin{lem}
  \label{l:singLocus} 
The singular locus of $\reduced{E}_{I}$ is contained in
$Z_{\reduced{E}_{I}}$. (Recall, that for an effective divisor
$D=\sum_\alpha \mu_\alpha C_\alpha$, where $C_\alpha$'s are reduced
irreducible curves, $\reduced{D}$ denotes the reduced associated
divisor $\sum_\alpha C_\alpha$.)
\end{lem}
 
\proof
Taking the restriction of the diagram \eqref{d:exteriorProduct} to
$\reduced{E}_{I}$ implies that, outside $Z_{\reduced{E}_{I}}$, the
sheaf $\Theta_{\reduced{E}_{I}}$ is isomorphic to the locally free
sheaf $\Oo_{\reduced{E}_{I}}(L)$.  By a result of Lipman, \cite{L},
this implies that $E^{red}_{I}$ is smooth outside the subscheme
$Z_{E^{red}_{I}}$.
\qed

Next we turn to the Zariski decomposition of $E$.  We first determine
the curves intersecting $E$ negatively.

\begin{lem}
  \label{l:nefForE}
Assume $K_X$ to be ample.  If $C$ is a reduced irreducible curve on
$X$ such that $C\cdot E<0$, then $C$ is a smooth rational component of
$E_{I}$ satisfying $C^2=-3$, $C\cdot E=-1$, and $C\cdot L=2$.
Furthermore, $C$ does not intersect any other component of $E_{I}$.
\end{lem}

\proof
The inequality $C\cdot E<0$ implies that $C$ is a component of $E$ and
that $C^2<0$.  Since $K_X$ is ample,
\begin{equation}
  \label{eq:C.L}
  C\cdot L = C\cdot (K_X-E) \geq 2,
\end{equation}
so by \eqref{eq:typeII}, $C$ must be a component of \typeI\ and
according to \eqref{eq:typeI}, $C\cdot L\leq2$.  This inequality
and \eqref{eq:C.L} yield $C\cdot L=2$ and $Z_C=\emptyset$.
The latter, combined with Lemma~\ref{l:singLocus}, implies that $C$ is
a smooth rational curve with $C\cdot K_X =1$, \ie $C$ is a
$(-3)$-curve, and that $C$ intersects no other component of $E_I$.
\qed

\begin{rem}
  \label{Knonample}
If we drop the assumption of ampleness for $K_X$, then the inequality
in \eqref{eq:C.L} implies that the intersection $C\cdot E$ is either
$-2$ or $-1$.  The first case implies that $C$ is a $(-2)$-curve, while
the second leads to two possibilities: either $C\cdot L=1$ and $C$ is
again a $(-2)$-curve, or $C\cdot L=2$ and $C$ is a $(-3)$-curve. 
\end{rem}

\begin{pro}
  \label{p:ZariskiDecompositionForE}
Let $K_X$ be ample and set $\RRR_{-1}=\{C \mid C\cdot E=-1\}$.  Then
the Zariski decomposition of $E$ is $E=E^++E^-$, where
\[
  E^- = \frac{1}{3}\,\sum_{C\in\RRR_{-1}}C.
\]
In particular, $E^+\neq0$.  
\end{pro}

\proof
We know that $E\neq0$.  
Let $E^+=E-E^-$, with $E^-$ defined as above.  It is sufficient to
show that $C\cdot E^+=0$, for any $C\in\RRR_{-1}$, and that $E^+$ is
nef.  The non-triviality of $E^+$ then follows immediately.

We begin by checking that the curves in $\RRR_{-1}$ are orthogonal to
$E^+$.  To do this, we observe that if $\Gamma\in\RRR_{-1}$, then
$\Gamma\cdot L=2$; from \eqref{eq:typeI}, it follows that the
subscheme $Z_\Gamma$ is empty.  Using Lemma~\ref{l:singLocus}, we see
that $\Gamma$ intersects no other curve from $\RRR_{-1}$, hence
\[
  \Gamma\cdot E^+ = 
  \Gamma\cdot \bigg(E-\frac{1}{3}\,\sum_{C\in\RRR_{-1}}C\bigg)
  =\Gamma\cdot E - \frac{1}{3}\,\Gamma^2
  =0.
\]

Next we show that $E^+$ is nef.  Let us suppose that there exists a
reduced irreducible curve $\Gamma$ such that $\Gamma\cdot E^+<0$.
By the above argument, $\Gamma\notin\RRR_{-1}$ and
\[
  0 \leq \Gamma\cdot E 
  = \Gamma\cdot E^+ + \Gamma\cdot E^- < \Gamma\cdot E^-.
\]
Putting together this inequality and the last assertion of
Lemma~\ref{l:nefForE}, we deduce that $\Gamma$ is a \typeII\ component
of $E$.  The rest of the argument is an estimate of $\Gamma\cdot E$.
To begin with,
\[
  \Gamma\cdot E = \Gamma\cdot(E_I+E_{II}) 
  \geq \Gamma\cdot E_I+m_\Gamma\Gamma^2,
\]
where $m_\Gamma$ is the multiplicity of $\Gamma$ in $E$.  Using
\eqref{eq:typeII}, 
\[
  \Gamma\cdot E 
  \geq \Gamma\cdot E_I+m_\Gamma(\Gamma\cdot L+\deg Z_\Gamma)
  \geq \Gamma\cdot E_I+m_\Gamma \Gamma\cdot K_X-m_\Gamma\Gamma\cdot E,
\]
hence
\[
  \Gamma\cdot E
  \geq \frac{m_\Gamma}{m_\Gamma+1}\,\Gamma\cdot K_X
  +\frac{1}{m_\Gamma+1}\,\Gamma\cdot E_I
  \geq \frac{m_\Gamma}{m_\Gamma+1}\,\Gamma\cdot K_X
  +\frac{1}{m_\Gamma+1}\sum_{C\in\RRR_{-1}}m_C\,\Gamma\cdot C.
\]
Since by Lemma~\ref{l:nefForE}, every $C\in\RRR_{-1}$
intersects no other curve in $E_I$,
\[
  -1 = C\cdot E 
  = -3m_C+C\cdot E_{II} \geq -3m_C+m_\Gamma\,C\cdot\Gamma.
\]
Thus
\[
  m_C \geq \frac{1}{3}\,\bigg(m_\Gamma\,\Gamma\cdot C+1\bigg).
\]
Using this inequality and the last estimate for $\Gamma\cdot E$, we
arrive at
\[
\begin{aligned}
  \Gamma\cdot E
  &\geq \frac{m_\Gamma}{m_\Gamma+1}\,\Gamma\cdot K_X
  +\frac{1}{3}\,\frac{1}{m_\Gamma+1}
  \sum_{C\in\RRR_{-1}}(m_\Gamma\,\Gamma\cdot C+1)\,\Gamma\cdot C \\
  &\geq \frac{m_\Gamma}{m_\Gamma+1}\,\Gamma\cdot K_X
  +\frac{1}{3}\,\sum_{C\in\RRR_{-1}}\Gamma\cdot C \\
  &= \frac{m_\Gamma}{m_\Gamma+1}\,\Gamma\cdot K_X
  + \Gamma\cdot E^-,
\end{aligned}
\]
which leads to the contradiction
\[
  0 > 
  \Gamma\cdot E^+ 
  \geq \frac{m_\Gamma}{m_\Gamma+1}\,\Gamma\cdot K_X > 0.
\]
\qed

\begin{rem*}
The preceding argument shows that every \typeII\ component $\Gamma$ of $E$
intersects $E^+$ positively.  More precisely, 
\[
\Gamma\cdot E^+ \geq \frac{m_\Gamma}{m_\Gamma+1}\,\Gamma\cdot K_X
\]
where $m_\Gamma$ is the multiplicity of $\Gamma$ in $E$.
\end{rem*}

Putting together Lemma~\ref{l:typeIIandLminus} with the property
\eqref{eq:typeI} of curves of \typeI, we obtain the
following.

\begin{pro}
  \label{p:mainOnE}
Assume $\alpha_X < \dfrac{5}{6}$ and 
let $C$ be an irreducible component of $E$.  Then
\listspace
\begin{itemize}
\item 
either $C\cdot L\leq0$ and then $C$ is orthogonal to $L^+$, the
positive part of the Zariski decomposition of $L$,
\item
or $C\cdot L>0$ and then $C$ is a smooth rational curve of \typeI\
with $C\cdot L =1$ or $2$.
\end{itemize}
\end{pro}

\proof
We consider $C$ according to the sign of its intersection product with
$L$.  If $C\cdot L\leq0$, we use the Zariski decomposition of 
$L = L^+ + L^-$ in Lemma \ref{l:L2-LK} to deduce that either $C$ is a
component of $L^-$, in which case $C\cdot L^+=0$, or it is not, in
which case
\[
  0 \leq C\cdot L^+ 
  = C\cdot(L-L^-) \leq -C\cdot L^- \leq 0,
\]
hence again $C\cdot L^+=0$.  

If $C\cdot L>0$, then, by Lemma~\ref{l:typeIIandLminus}, the curve $C$
is of \typeI, and by \eqref{eq:typeI},
\[
  0 \leq 2g(\Ctilde) \leq
  2-C\cdot L -\deg Z_C < 2.
\]
This implies that $C\cdot L\leq2$ and $g(\Ctilde)=0$.  Furthermore, by
Lemma \ref{l:singLocus}, $C$ is smooth outside the
subscheme $Z_C$, while the inequality \eqref{eq:typeI},
\[
  \deg (Z_C) \leq \deg\eta_C^\ast(Z_C) \leq 1,
\]
implies that $Z_C$ is either empty or a smooth point of $C$.
\qed

\begin{cor}
  \label{c:CAndLPlus}
Assume $K_X$ ample and $\alpha_X < \dfrac{5}{6}$. Then for every
irreducible component $C$ of $E$,
\[
  C\cdot L^+\leq 1
\]
Furthermore, if $C\cdot L^+>0$, then $C$ is a smooth rational curve of
type I such that
\[
  C\cdot L =j 
  \textq{and}
  C\cdot\reduced{E}_{II}=j-1,
  \quad j\in\{1,2\}.
\]
\end{cor}

\proof
Let $C$ be an irreducible component of $E$.  Using the previous
proposition, we know that $C\cdot L^+=0$ unless $C$ is a smooth
rational curve of \typeI\ such that $C\cdot L=j$, $j\in\{1,2\}$.  Let
$C$ be such a rational curve with $C\cdot L^+>0$.  Clearly $C$ is not
a component of $L^-$, hence $C\cdot L^-\geq0$, and we have
\begin{equation}
  \label{eq:C.L+}
  C \cdot L^+ = j - C\cdot L^-.
\end{equation}
In view of Lemma~\ref{l:typeIIandLminus}, if $j=1$, the curve $C$
cannot intersect any \typeII\ component of $E$, \ie
$C\cdot\reduced{E}_{II} =0$.  

We now turn to the case $j=2$.  First we
will show that $C\cdot E_{II} > 0$.  To see this, notice that
\begin{equation}
  \label{C.EII}
  C\cdot E_{II} = C\cdot E - C\cdot E_{I} \geq -1 - C\cdot E_{I},
\end{equation}
where the last inequality uses Lemma~\ref{l:nefForE}.
From $C\cdot L=2$ and the inequality \eqref{eq:typeI}, it follows that
$Z_C=\emptyset$.  Hence, according to Lemma~\ref{l:singLocus}, $C$
intersects no other component of $E_{I}$, \ie  
$C\cdot E_{I} = m_C C^2$, where $m_C$ is the multiplicity of $C$ in
$E_I$.  Substituting it in \eqref{C.EII} leads to
\[
  C\cdot E_{II} = -1 - m_C C^2 > 0.
\]
Thus, there are irreducible components $\Gamma$ in $E_{II}$
intersecting $C$ non-trivially.  Since, by
Lemma~\ref{l:typeIIandLminus}, each such $\Gamma$ appears in $L^-$
with multiplicity $\mu_{\Gamma}\geq1$, the equality \eqref{eq:C.L+}
implies that there exists a unique $\Gamma$ in $E_{II}$ intersecting
$C$ and $\Gamma\cdot C = 1$.  Hence $C\cdot\reduced{E}_{II} =1$ and 
$C\cdot L^+\leq 1$.
\qed

\begin{remark}
  \label{r:obst}
From Corollary \ref{c:CAndLPlus} it follows that
the main reason for the non-vanishing of $H^1(\Theta_X (-K_X))$
resides in the existence of curves $C$ of \typeI\ with $C\cdot L^+>0$.
Indeed, suppose that there would be no such curve.  Then 
$C\cdot L^+=0$ for every component $C$ of $E$.  This gives
$E^+\cdot L^+=0$.  Since $E^+$ is nef and non-zero, this is
impossible.
\end{remark}

\section{The proof of Theorem~\ref{t:theResult}}

We now assume that $K_X$ is ample and $\alpha_X<5/6$.  We seek to
obtain a contradiction to the decomposition
\[
K_X =L + E
\]
of Lemma~\ref{l:decomp}.  Our main idea is to compare $E^2$ and 
$E\cdot L^+$.  This will be realized by using the properties of the
divisors $L$ and $E$ accumulated in the previous sections.  In view of
Remark~\ref{r:obst}, we need to deal with the smooth rational curves
$C$ of \typeI\ subject to $C\cdot L^+ >0$.  By
Corollary~\ref{c:CAndLPlus}, these are components of $E_I$ satisfying
either $C\cdot L =1$, or $C\cdot L =2$.  For this reason, we introduce
the sets
\begin{equation}
  \label{setRj}
  \RRR^j=\{C \mid C \subset E^+, C \cdot L=j\},
  \qquad j=1,2.
\end{equation}
We have seen in Proposition~\ref{p:mainOnE} that every curve in
$\RRR^j$ is smooth and rational.  For any integer $k$, set
\[
  \RRR_k^j = \{C \in \RRR^j \mid C \cdot E = k\}.
\]
Since every $C\in\RRR^j$ verifies $0<C\cdot K_X=j+C\cdot E$, the sets
$\RRR^j$ are partitioned by the subsets $\RRR_k^j$ as follows:
\[
  \RRR^1 = \RRR_0^1 \cup \RRR_{\geq1}^1
  \textq{and}
  \RRR^2 = \RRR_{-1}^2 \cup \RRR_0^2 \cup \RRR_{\geq1}^2
\]
where $\displaystyle{\RRR_{\geq1}^j =\bigcup_{k\geq1} \RRR_k^j}$.

\begin{rem*}
$\RRR_{-1}^2=\RRR_{-1}$, \cf Lemma~\ref{l:nefForE}.
\end{rem*}

Using Proposition~\ref{p:mainOnE} and the above partitions of the sets
$\RRR_j$,
\begin{equation}
  \label{E-decomp}
  E = Q + \sum_{C\in\RRR_0^1}m_C C
  + \sum_{C\in\RRR_0^2}m_C C + \sum_{C\in\RRR_{-1}^2}m_C C,
\end{equation}
where $Q$ contains all the components lying in $\RRR_{\geq1}^j$
($j=1,2$) and those for which $C\cdot L^+=0$.

\begin{lem}
  \label{l:noMinus3}
If $X$ contains no $(-3)$-curve, then $\RRR_0^1=\RRR_{-1}^2=\emptyset$
and $E$ is nef.
\end{lem}

\proof
If $C$ is a smooth rational curve in $\RRR_0^1\cup \RRR_{-1}^2$, then
$C\cdot K_X=1$ and $C^2=-3$. Hence the first assertion. The second one 
follows immediately
from the Zariski decomposition of $E$ in Proposition 
\ref{p:ZariskiDecompositionForE}.
\qed

From now on we also assume that $X$ contains no $(-3)$-curves.  With
this hypothesis, our idea of comparison of $E^2$ and $E\cdot L^+$
takes the following form:

\begin{lem}
  \label{l:theInequality}
$E^2 \geq E\cdot L^+$.
\end{lem}

Let us assume the lemma and complete the argument.  First, notice that
this inequality implies $E^2>0$.  Then, by the Hodge index theorem,
\[
  E^2 \geq (L^+)^2 \geq L^2 \geq \frac{3}{2}\,(1-\alpha_X)\,K_X^2,
\] 
where the last inequality comes from \eqref{eq:LSquareAndKSquare}.
Arguing as in the proof of Lemma~\ref{eq:mainInequality}, we obtain
\[
  E\cdot K_X \geq \sqrt{\frac{3}{2}\,(1-\alpha_X)}\,K_X^2
  \textq{and}
  L^+\cdot K_X \geq \sqrt{\frac{3}{2}\,(1-\alpha_X)}\,K_X^2.
\]
Thus
\[
  K_X^2 \geq (L^++E)\cdot K_X \geq \sqrt{6(1-\alpha_X)}\,K_X^2,
\]
yielding $\alpha_X\geq 5/6$, which contradicts the
hypotheses.  At this stage Theorem \ref{t:theResult} is proved and we
turn to

\likeproof[Proof of Lemma~\ref{l:theInequality}]
Combining \eqref{E-decomp}, Lemma~\ref{l:noMinus3}, and
Corollary~\ref{c:CAndLPlus}, we obtain
\begin{equation}
  \label{eq:E-Q}
  (E-Q)\cdot L^+ 
  = \bigg(\sum_{C\in\RRR_0^2}m_CC\bigg)\cdot L^+
  \leq \sum_{C\in\RRR_0}m_C,
\end{equation}
where 
\[
  \RRR_0 = \{C\in\RRR_0^2 \mid C\cdot L^+ >0\}.
\]
It is clear that the curves in $\RRR_0$ do not lie in $L^-$.
Furthermore, by Corollary~\ref{c:CAndLPlus}, every $C$ in $ \RRR_0$
intersects a unique irreducible component of \typeII.  Denote this
component by $\Gamma_C$.  The correspondence $C\mapsto \Gamma_C$ gives
rise to the map
\[
  \rho : \RRR_0 \lra \Irr(E_{II}),
\]
where $\Irr(E_{II})$ stands for the set of irreducible components of
$E_{II}$.  Using $\rho$, we will be able to replace the multiplicities
$m_C$ in \eqref{eq:E-Q} by the ones of the curves
$\Gamma\in\Irr(E_{II})$.  This is done as follows.

Let $\Gamma$ be a curve in $\Irr(E_{II})$ and set
$\RRR_0(\Gamma)=\rho\inv(\Gamma)$.  Then for every $C\in\RRR_0(\Gamma)$
we have
\[
  0 = C\cdot E = C\cdot(E_I+E_{II}) = -4m_C+m_\Gamma,
\]
since $C^2 =-4$, $C\cdot\Gamma=1$, and $C$ intersects no other
component of $E$.  Thus $m_\Gamma=4m_C$.  Substituting it in
\eqref{eq:E-Q}, yields
\begin{equation}
  \label{eq:final}
  (E-Q)\cdot L^+ \leq  \sum_{C\in\RRR_0}m_C 
  = \frac{1}{4}\sum_{\Gamma\in\Irr(E_{II})} n_\Gamma\,m_\Gamma,
\end{equation}
where $n_\Gamma=\card(\rho\inv(\Gamma))$.

\paragraph{Claim}
If\/ $\Gamma$ is in $\Irr(E_{II})$, then 
\[
  \Gamma\cdot E\geq\dfrac{1}{4}\,n_\Gamma.
\]

Assuming the claim, let us complete the proof of
Lemma~\ref{l:theInequality}.  We combine the inequality of
the claim and \eqref{eq:final} to obtain
\[
  (E-Q)\cdot L^+ \leq 
  \frac{1}{4}\sum_{\Gamma\in\Irr(E_{II})} n_\Gamma\,m_\Gamma
  \leq \sum_{\Gamma\in\Irr(E_{II})}m_\Gamma \Gamma\cdot E
  = E_{II}\cdot E.
\]
This and the definition of $Q$ yield 
\[
  E\cdot L^+ 
  \leq Q\cdot L^+ + E_{II}\cdot E 
  = \bigg(
  \sum_{C\in\RRR_{\geq1}^1}m_CC+\sum_{C\in\RRR_{\geq1}^2}m_CC
  \bigg)\cdot L^+ + E_{II}\cdot E.
\]
Using Corollary~\ref{c:CAndLPlus} and the fact that $E$ is nef 
(Lemma \ref{l:noMinus3}), we
obtain
\[
\begin{split}
  E\cdot L^+
  &\leq \sum_{C\in\RRR_{\geq1}^1\cup\RRR_{\geq1}^2}m_C+ E_{II}\cdot E\\
  &\leq \sum_{C\in\RRR_{\geq1}^1\cup\RRR_{\geq1}^2}m_CC\cdot E
  + E_{II}\cdot E\\
  &\leq E_I\cdot E + E_{II}\cdot E\\
  &= E^2.
\end{split}
\]
\qed

\likeproof[Proof of the claim]
We consider the divisor
$\Bigg(\Gamma+\dfrac{1}{4}\sum_{C\in\RRR_0(\Gamma)}C\Bigg)$ and argue
according to its self-intersection,
\[
  \Bigg(\Gamma+\dfrac{1}{4}\sum_{C\in\RRR_0(\Gamma)}C\Bigg)^2
  = \Gamma^2+\frac{n_\Gamma}{4}.
\]
 If $\Gamma^2\leq-\dfrac{n_\Gamma}{4}$, then we obtain
\[
-\frac{n_\Gamma}{4} \geq \Gamma^2 \geq  \Gamma\cdot L=\Gamma\cdot(K_X-E)
\]
where the second inequality uses the fact that $\Gamma$ is of \typeII\
and hence subject to the inequality \eqref{eq:typeII}.  Thus we obtain
\[
  \Gamma\cdot E \geq \Gamma\cdot K_X-\Gamma^2 \geq 
  \dfrac{n_\Gamma}{4} + \Gamma\cdot K_X > \dfrac{n_\Gamma}{4}.
\]
If $\Gamma^2 > -\dfrac{n_\Gamma}{4}$, we use the fact that the divisor
$\Bigg(\Gamma+\dfrac{1}{4}\sum_{C\in\RRR_0(\Gamma)}C\Bigg)$ enters $E$
with multiplicity $m_{\Gamma}$, \ie
\[
  E = 
  m_{\Gamma}\Bigg(
  \Gamma+\dfrac{1}{4}\sum_{C\in\RRR_0(\Gamma)}C
  \Bigg) + R_\Gamma,
\]
where $R_\Gamma$ is the residual part of $E$.  This implies
\begin{equation}
  \label{eq:Ga.E}
  \begin{aligned}
    \Gamma\cdot E
    &= m_{\Gamma} \Gamma\cdot
    \Bigg(\Gamma+\dfrac{1}{4}\sum_{C\in\RRR_0(\Gamma)}C\Bigg)
    + \Gamma\cdot R_\Gamma\\ 
    &\geq m_\Gamma\,\Gamma^2+\frac{n_\Gamma m_\Gamma}{4}
  \end{aligned}
\end{equation}
By \eqref{eq:typeII}, 
\[
  \Gamma^2 \geq \Gamma\cdot L = \Gamma\cdot K_X-\Gamma\cdot E.
\]
Substituting into \eqref{eq:Ga.E},
\[
  \Gamma\cdot E \geq  m_\Gamma\,\Gamma^2
    +\frac{n_\Gamma m_\Gamma}{4}
  \geq m_\Gamma(\Gamma\cdot K_X-\Gamma\cdot E)
  +\frac{n_\Gamma m_\Gamma}{4}.
\]
Solving for $\Gamma\cdot E$, we obtain 
\[
  (m_\Gamma+1)\,\Gamma\cdot E \geq m_\Gamma\,\Gamma\cdot K_X
  +\frac{n_\Gamma m_\Gamma}{4}.
\]
Summing up with \eqref{eq:Ga.E} leads to
\begin{equation}
  \label{eq:Ga.E-1}
  (m_\Gamma+2)\,\Gamma\cdot E 
  \geq m_\Gamma(\Gamma\cdot K_X+\Gamma^2)
  +\frac{n_\Gamma m_\Gamma}{2}
\end{equation}
and we continue to argue according to the sign of the expression 
$(\Gamma\cdot K_X+\Gamma^2)$.

\case
If $\Gamma\cdot K_X+\Gamma^2 \geq 0$, then the inequality
\eqref{eq:Ga.E-1} yields
\[
  \Gamma\cdot E \geq \frac{n_\Gamma m_\Gamma}{2(m_\Gamma+2)}
  \geq \frac{n_\Gamma}{3}
  > \frac{n_\Gamma}{4},
\]
where the second inequality follows from $m_\Gamma=4m_C\geq 4$.  

\case
If $\Gamma\cdot K_X+\Gamma^2 < 0$, then 
$\Gamma\cdot K_X+\Gamma^2 = -2$ and $\Gamma$ is a smooth rational
curve.  Substituting this relation in \eqref{eq:Ga.E-1} yields
\[
  (m_\Gamma+2)\,\Gamma\cdot E 
  \geq -2m_\Gamma+\frac{n_\Gamma m_\Gamma}{2}.
\]
So
\begin{equation}
  \label{eq:Ga.E-2}
  \Gamma\cdot E
  \geq \frac{m_\Gamma}{2(m_\Gamma + 2)}\,(n_\Gamma -4).
\end{equation}   
Since $\Gamma$ is a smooth rational curve, $\Gamma^2 \leq -4$---the
surface $X$ contains neither $(-2)$-curves nor $(-3)$-curves.  Together 
with the assumption $\Gamma^2>-\dfrac{n_\Gamma}{4}$, this inequality
implies that $n_\Gamma>16$.  Using \eqref{eq:Ga.E-2} and
$m_\Gamma\geq4$, we obtain
\[
  \Gamma\cdot E 
  \geq \frac{1}{3}\,(n_\Gamma -4)
  \geq 
  \frac{n_\Gamma}{4}+\frac{n_\Gamma}{12}-\frac{4}{3} > \frac{n_\Gamma}{4}.
\]
\qed

It is reasonable to guess that once the group $H^1(\Theta_X(-m_0K_X))$
vanishes for some $m_0$, the same should hold for all $m>m_0$.  This
is justified by the following:

\begin{pro}
  \label{p:for-mK}
Let $X$ be a smooth minimal complex surface of general type with
$p_g(X)\geq1$ and 
$\alpha_X<1$.  Let $c$ be a global section of $\Oo_X(K_X)$
and let $C=\{c=0\}$ be the corresponding divisor.  If $C$ is reduced
and irreducible, then the homomorphism 
\[
  H^1(\Theta_X(-(m+1)K_X)) \stackrel{c}{\lra} 
  H^1(\Theta_X(-mK_X))
\]
of the multiplication by $c$ is injective for all $m\geq0$. 
\end{pro}

\proof
Let $\xi$ be a non-zero class in $H^1(\Theta_X(-(m+1)K_X))$ such hat
$c\cdot \xi=0$.  Using the identification 
\[
  H^1(\Theta_X(-(m+1)K_X)) \cong 
  \Ext^1(\Omega_X(mK_X),\Oo_X(-K_X)),
\]
we view $\xi$ as the corresponding non-trivial extension class, \ie
the short exact sequence of vector bundles on $X$
\[
  0 \lra 
  \Oo_X(-K_X) \lra
  \Tt_\xi \lra
  \Omega_X(mK_X)) \lra 0.
\]
We dualize and consider the long exact sequence of cohomology groups
\[
  0 \lra H^0(\Tt_\xi^\ast)
  \lra H^0(\Oo_X(K_X))
  \stackrel{\xi}{\lra} H^1(\Theta_X(-mK_X)) \lra\cdots
\]
where the coboundary map is given by the cup-product with the extension
class $\xi$.  The assumption $c\cdot\xi=0$ implies that the section $c$
is the image of a non-zero section $s\in H^0(\Tt_\xi^\ast)$.  This
gives rise to the following commutative diagram.
\begin{equation}
  \label{d:dual}
  \UseTips
  \xymatrix@W=5mm@H=5mm{
    && 0 \ar[d] \\
    && \Oo_X \ar[d]_s\ar[rd]^{c} \\
    0  \ar[r] & \Theta_X(-mK_X) \ar[r]&
    \Tt_\xi^\ast \ar[r]\ar[d] & \Oo_X(K_X)\ar[r] & 0 \\
    && \coker(s) \ar[d] \\
    && 0
  }
\end{equation}
Furthermore, since the divisor $C=\{c=0\}$ is reduced
and irreducible, and the extension $\xi$ is non-trivial, it follows
that the sheaf $\coker(s)$ in \eqref{d:dual} is torsion free.
Dualizing again, we obtain the commutative diagram 
\begin{equation}
  \label{d:for-mK}
  \UseTips
  \xymatrix@W=5mm@H=5mm{
    && 0 \ar[d] \\
    && \Ff \ar[d]\ar[rd]^{\phi_\xi} \\
    0 \ar[r] & \Oo_X(-K_X) \ar[r] \ar[rd]_{\psi_\xi}&
    \Tt_\xi \ar[r]\ar[d] & \Omega_X(mK_X) \ar[r] & 0 \\
    && \Jj_Z \ar[d] \\
    && 0
  }
\end{equation}
where $Z$ is the subscheme of zeros of $s$, $\Jj_Z$ is its sheaf of
ideals, and $\Ff=(\coker(s))^\ast$ is a locally free sheaf of rank
$2$. 

From \eqref{d:for-mK} we deduce the Chern invariants of $\Ff$:
\[
  c_1(\Ff)=2mK_X
  \textq{and}
  c_2(\Ff) = c_2(X)+(m^2-m-1)K_X^2-\deg(Z).
\]
This and the assumption $\alpha_X<1$ imply
\[
  c_1^2(\Ff)-4c_2(\Ff) = 
  4\big((m+1)K_X^2-c_2(X)\big)+4\deg(Z) > 0.
\]
Thus $\Ff$ is Bogomolov unstable and its Bogomolov destabilizing
subsheaf $\Oo_X(F)$ gives rise to the exact sequence 
\begin{equation}
  \label{eq:destabilizingSequence}
  0 \lra
  \Oo_X(F) \lra
  \Ff \lra
  \Jj_A(2mK_X-F) \lra 0,
\end{equation}
where $A$ is a $0$-dimensional subscheme and $\Jj_A$ its sheaf of
ideals.  The fact that $\Oo_X(F)$ is the Bogomolov destabilizing sheaf
also implies that the divisor
\[
  c_1(\Oo_X(F))-\frac{1}{2}\,c_1(\Ff) = F-mK_X
\]
is in the positive cone of $NS(X)_\QQ$.  But the morphism $\phi_\xi$ in
\eqref{d:for-mK} combined with \eqref{eq:destabilizingSequence}
induces the monomorphism 
\[
  \Oo_X(F-mK_X) \lra \Omega_X.
\]
Hence $\Omega_X$ has a subsheaf of rank $1$ and Iitaka dimension $2$
which is impossible. 
\qed

\begin{cor}
  \label{c:vanish}
Let $X$ be a smooth minimal surface of general type such that
$\alpha_X=c_2(X)/K_X^2<5/6$ and the canonical linear system $|K_X|$
contains a reduced irreducible member.  If $X$ contains no smooth
rational curves with self-intersection $-2$ or $-3$, then
\[
  H^1(\Theta_X(-mK_X)) = 0, \textq{for all} m\geq1.
\] 
\end{cor}

\proof
Follows immediately from Theorem~\ref{t:theResult} and
Proposition~\ref{p:for-mK}.
\qed

\section{ The irregularity of surfaces in $\PP^4$}

The aim of this section is to study the geography of
surfaces of general type $X$ in $\PP^4$ with fixed numerical
invariants, $d=\deg X$, $K_X^2$, and $\chi(\Oo_X)$.  To this end
we consider the set
\begin{equation}
  \label{eq:Sdn}
  \sS_d(\chi,\beta) := \{ X \mid
  \begin{aligned}[t]
    X\subset\PP^4\, 
    &\text{ smooth surface of general type},\\
    &\deg X=d,\, \chi(\Oo_X)=\chi,\,
    K_X^2=\beta \chi(\Oo_X)\},
  \end{aligned}
\end{equation}
where $\beta\leq9$. There is a general result due to Decker and
Schreyer, \cite[Proposition 3.11]{DeSch} (see also \cite[Proposition
3]{E-P}), which establishes an upper bound on $d$ for each given value
of $\chi$.  Here we will explore the irregularity of surfaces in
$\sS_d(\chi,\beta)$.  More precisely, for the surfaces in the
subset $\sS_d^\circ(\chi,\beta)\subset\sS_d(\chi,\beta)$ defined by
\[
  \sS_d^\circ(\chi,\beta) := 
  \{ X\in\sS_d(\chi,\beta) \mid 
  X\text{ minimal and with no irrational pencil} \}
\]
we are able to show the following:

\begin{thm}
  \label{bound-q}
Given a positive integer $\chi$ and a rational $\beta$,
$2\leq \beta\leq9$, there exists a degree $d(\chi,\beta)$ such that,
if $d>d(\chi,\beta)$, then every surface
$X\in\sS_d^\circ(\chi,\beta)$ has the irregularity at most $3$.
\end{thm}

If the irregularity is equal to $3$, our argument (see
Lemma~\ref{l:4-5}) shows that $X$ must be contained in a quintic
$3$-fold in $\PP^4$ and must satisfy one of the following conditions:

\begin{enumerate}\it
\item[(i)]
  $\Omega_X$ is generated by global sections, and
  $K^2_X=6\chi(\Oo_X)$;
\item[(ii)]
  the canonical divisor has the form
\[
  K_X = \mu H + F + \sum_{C\cdot H=2}m_C C + \sum_{l\cdot H=1} m_l l,
\]
where $\mu\geq 1$, $F$ is an effective divisor, and the sums run over
collections of irreducible conics ($C\cdot H=2$) and lines 
($l\cdot H=1$) contained in $X$.
\end{enumerate}

A proof of this theorem is a consequence of several lemmas.  Though
the proofs of some of those results are somewhat lengthy and sometimes
involve considerations of several cases, the main principle behind all
the arguments is quite transparent: the exploiting of the {\it
non-vanishing} of $H^1(X,\Theta_X(-K_X))$.  This non-vanishing gives
rise, as in Lemma \ref{l:decomp}, to a decomposition of the canonical
divisor $K_X$.  Then we seek to relate this
decomposition to the linear system embedding $X$ in $\PP^4$.

The non-vanishing of $H^1(X,\Theta_X(-K_X))$ turns out to be a
cohomological version of the following result of Ellingsrud and
Peskine in \cite{E-P}.

\begin{lem}
  \label{l:EP}
Given the integer
$\chi$ and the rational $\beta$, $\beta\leq9$, there exists a degree
$d(\chi,\beta)$ such that for all $d>d(\chi,\beta)$ and for every
surface $X$ in $\sS_d(\chi,\beta)$, there is a $3$-fold in $\PP^4$ of
degree $\leq 5$ containing $X$.
\end{lem}

\proof
If $X$ is embedded in $\PP^4$ by the line bundle $\Oo_X(H)$, the
numerical invariants of $X$ satisfies the well-known formula (see \eg
\cite{H}, p.434),
\begin{equation}
  \label{eq:d-and-Chern}
  d^2 -10d -5H\cdot K_X = 2K^2_X - 12\chi(\Oo_X), 
\end{equation}
or, using the geometric genus of a hyperplane section, $g_H$, the
adjunction formula, and the notation above,
\begin{equation}
  \label{d-and-Chern}
  d^2-5d-10(g_H-1)-2(\beta-6)\chi = 0.
\end{equation}

By a result of Gruson and Peskine in \cite{G-P} (the Halphen's
bound), if $d\geq 30$ and a hyperplane section $H$ {\it is not}
contained in a surface of degree $5$, then
\[
  g_H-1 \leq \frac{d^2 +12d}{12}.
\]
Putting this together with \eqref{d-and-Chern} yields 
\[
  d^2-5d-2(\beta-6)\chi \leq \frac{5}{6}\,(d^2 +12d),
\]
or, equivalently,
\[
  d^2-90d-12(\beta-6)\chi \leq 0.
\]
Set $d(\chi,\beta)=29$ if the discriminant of the left hand side is
negative. Otherwise, we solve the inequality for $d$ and set
\[
  d(\chi,\beta) =  \max(29,45+\sqrt{45^2+12(\beta-6)\chi}).
\]
Then, if $d > d(\chi,\beta)$, the hyperplane sections of $X$ must be
contained in a surface of degree $\leq 5$.  This and the classical
result of Roth (\cite{Ro}, p.152, see also \cite{E-P}, (C), p.2) imply
that $X$ is contained in a $3$-fold of degree $\leq 5$ in $\PP^4$,
provided $d>d(\chi,\beta)$.
\qed

Let $\Jj_X$ be the ideal sheaf of $X$ in $\PP^4$.  Set
\begin{equation}
  \label{mX}
  m_X := \min\{m\in {\bf Z} \mid h^0 (\Jj_X (m)) \neq 0\},
\end{equation}
\ie $m_X$ is the minimal degree of hypersurfaces in $\PP^4$ containing
$X$. 

Let $\Nn_X$ be the normal bundle of $X$ in $\PP^4$ and let
$\Nn^{\ast}_X$ be its dual, the conormal bundle.  The normal bundle is
a rank $2$ bundle on $X$ with determinant bundle
$\det(\Nn_X)=\wedge^2\Nn_X=\Oo_X(K_X+5H)$, where $H$ is a hyperplane
section of $X$, \ie $\Oo_X(H)=\Oo_{\PP^4}(1)\otimes\Oo_X$.  This
implies
\begin{equation}
  \label{nor-conor}
  \Nn^\ast_X \cong \det(\Nn^\ast)\otimes\Nn_X = \Nn_X(-K_X-5H).
\end{equation}
By the definition of $m_X$ in \eqref{mX} and the identity 
$\Nn^\ast_X = \Jj_X / {\Jj}^2_X$, we deduce 
$h^0 (\Nn^\ast_X (m_X)) \neq 0$.  This together with the
identification in \eqref{nor-conor} yield
\[
  H^0(\Nn^\ast_X(m_X)) = H^0(\Nn_X(-K_X-(5-m_X)H)) \neq 0.
\]
This non-vanishing is related to the $H^1$ of the tangent bundle
$\Theta_X$ via the normal sequence of $X$ in $\PP^4$, 
\[
  0 \lra \Theta_X \lra \Theta_{\PP^4}\otimes\Oo_X \lra \Nn_X \lra 0.
\]
Tensoring with $\Oo_X(-K_X-(5-m_X)H)$ and passing to the the long
exact sequence of cohomology groups, we obtain 
\begin{equation}
  \label{coh-seq}
  H^0 (\Theta_{\PP^4}\otimes\Oo_X(-D)) \lra
  H^0 (\Nn_X (-D)) \stackrel{\delta_X}{\lra}
  H^1 (\Theta_X(-D))\,,
\end{equation}
where we set $D=K_X+(5-m_X)H$.  In view of Lemma \ref{l:EP}, we
are interested in the case $m_X\leq 5$.  The following result is a
cohomological interpretation of that lemma.

\begin{lem}
  \label{l:EP-coh}
If $m_X\leq5$, then the coboundary homomorphism $\delta_X$ in
\eqref{coh-seq} is injective unless $m_X =5$, $p_g(X)\leq 2$, and the
degree of $X$ in $\PP^4$ is at most $16$.
\end{lem}

\proof
The failure of injectivity of $\delta_X$ in \eqref{coh-seq} implies
the non-vanishing of 
\[
  H^0 (\Theta_{\PP^4}\otimes\Oo_X(-D))=
  H^0(\Theta_{\PP^4}\otimes\Oo_X(-K_X-(5-m_X)H)).
\]
Using the Euler sequence on $\PP^4$,
\[
  0 \lra \Oo_{\PP^4} \lra H^0(\Oo_{\PP^4}(1))^\ast\otimes\Oo_{\PP^4}(1)
  \lra \Theta_{\PP^4} \lra 0,
\]
it follows, by restriction to $X$, that either
$H^0(\Oo_X(-K_X-(4-m_X)H))$ or $H^1(\Oo_X(-K_X-(5-m_X)H))$ must be
nonzero.  But for $m_X\leq 4$ both groups vanish: the first one
because $X$ is of general type, and the second one because the divisor
$(K_X+(5-m_X)H)$ is nef and big.  Hence, we must have $m_X=5$, \ie
$H^0(\Theta_{\PP^4}\otimes\Oo_X(-K_X )\neq0$.

From the Euler sequence, this non-vanishing can occur either if
$H^0(\Oo_X(H-K_X))\neq 0$, or if
\begin{equation}
  \label{ker}
  \ker\Big(
  H^1(\Oo_X(-K_X)) \lra H^0(\Oo_X(H))^\ast \otimes H^1 (\Oo_X(H-K_X))
  \Big) \neq 0.
\end{equation}
We consider each of these two cases.

\case
If $H^0(\Oo_X(H-K_X))\neq 0$, then 
\begin{align}
  &\textq{either} H=K_X,\nonumber\\
  &\textq{or} H-K_X = \Gamma,
  \textq{with $\Gamma$ a non-zero effective
  divisor.} \label{p2}
\end{align}
We claim that the former is impossible.  Indeed, if $H=K_X$, then by
\cite{B-C}, $X$ is a complete intersection of two hypersurfaces
in $\PP^4$.  Since $m_X=5$, it follows that the degrees of these
hypersurfaces are $5$ and $n\geq 5$.  But then, by the adjunction
formula for $K_X$, we have $K_X = nH$, which contradicts
$H=K_X$.  Thus \eqref{p2} holds and we deduce
\begin{equation}
  \label{pg3}
  p_g(X) = h^0(\Oo_X(K_X)) = h^0(\Oo_X(H-\Gamma)) \leq 3.
\end{equation}
Furthermore, the inequality must be strict, \ie $p_g(X)\leq 2$.
Indeed, if the equality holds in \eqref{pg3}, then $\Gamma$ is a line
in $\PP^4$.  Intersecting both sides in \eqref{p2} with the line
$\Gamma$, we obtain
\[
  1=H\cdot\Gamma = (K_X+\Gamma)\cdot\Gamma = -2,
\]
which is absurd.  Thus $p_g(X)\leq 2$ and it remains to bound the
degree $d$ of $X$.  From \eqref{p2},
\[
  H\cdot K_X = H^2 - \Gamma\cdot H < d.
\]
Substituting this in the relation \eqref{d-and-Chern}, we obtain
\begin{equation}
  \label{d-Chern}
  d^2-15d < 2K^2-12\chi(\Oo_X).
\end{equation}
Since $p_g(X) \leq 2$, it follows $\chi(\Oo_X)\leq 3$.  This and the
Bogomolov-Miyaoka-Yau inequality imply 
\[
  2K^2 -12\chi(\Oo_X) \leq 18.
\]
Putting this together with \eqref{d-Chern}, we obtain 
$d\leq 16$.

\case
If $\ker\Big(
H^1(\Oo_X(-K_X)) \lra H^0(\Oo_X(H))^\ast \otimes H^1 (\Oo_X(H-K_X))
\Big) \neq 0$,
let $X_0$ be the minimal model of $X$.  Since $X$ is of general type,
$X_0$ is uniquely defined.  Let
\[
  f: X\lra X_0
\]
be the corresponding sequence of blowing-down maps.  The canonical
divisor of $X$ can be written as follows:
\begin{equation}
  \label{K-Kmin}
  K_X = f^\ast K_{X_0} + \Delta
\end{equation}
where $\Delta$ is the exceptional divisor of $f$.  In particular,
$\Delta$ is composed of rational curves on $X$ contracted to points by
$f$.  Thus $\Delta\cdot f^\ast K_{X_0}=0$ and $f^\ast K_{X_0}$ is nef
and big---$X_0$ is minimal and of general type.

The non-vanishing of $H^0(\Theta_{\PP^4}\otimes\Oo_X(-K_X ))$ implies
the non-vanishing of 
$H^0(\Theta_{\PP^4}\otimes \Oo_X(-f^\ast K_{X_0}))$.  From the Euler
sequence and the fact that $f^{\ast} K_{X_0} $ is nef and big, the
latter non-vanishing yields
\[
  H^0(\Oo_X(H-f^\ast K_{X_0})) \neq 0.
\]
As before, there are two possibilities, either $H=f^\ast K_{X_0}$, or
$H-f^\ast K_{X_0} =\Gamma$, with $\Gamma$ effective and nonzero.
Since $H$ is very ample, the first possibility yields $X_0 =X$ and
$H=K_X$ which implies the vanishing of $H^1 (\Oo_X (-K_X))$, contrary
to our assumption.  Hence the latter possibility must hold.

Arguing as in the first case, we deduce $p_g(X)\leq 3$.  Furthermore,
the inequality must be strict.  Indeed, if $p_g(X)=3$, then $\Gamma$
is a line.  Since $H-f^\ast K_{X_0} =\Gamma$, it follows that
\begin{equation}
  \label{H-K0line}
  H=K_X + \Gamma -\Delta.
\end{equation}
Intersecting the both sides with the line $\Gamma$ yields
\[
  1 = H\cdot\Gamma = (K_X+\Gamma)\cdot\Gamma -\Delta\cdot\Gamma 
  = -2-\Delta\cdot\Gamma.
\]
Hence $\Delta\cdot\Gamma = -3$ and the line $\Gamma$ is a component of
$\Delta$, \ie
\begin{equation}
  \label{Del-l}
  \text{the divisor $(\Delta-\Gamma)$ is effective}.
\end{equation} 
On the other hand, the condition \eqref{ker} implies that for every
$h\in H^0(\Oo_X(H))$, the cup-product
\[
  H^1(\Oo_X(-K_X)) \stackrel{h}{\lra} H^1( \Oo_X(H-K_X))
\]
has a nontrivial kernel.  This implies that 
$H^0(\Oo_C(H-K_X)) \neq 0$, for every $C$ in the linear system 
$|H|$.  By \eqref{H-K0line}, $H-K_X=\Gamma-\Delta$, hence
\begin{equation}
  \label{onC}
  H^0(\Oo_C(\Gamma-\Delta))\neq 0 
  \text{ for every $C\in|H|$}.
\end{equation}
This, together with \eqref{Del-l}, implies $\Delta=\Gamma$.
Substituting it in \eqref{H-K0line} leads to $H=K_X$ which is
impossible.  Thus $p_g(X)\leq 2$ and it remains to give the upper
bound on $d$.

From \eqref{onC}, it follows that $H\cdot(\Gamma -\Delta)\geq 0$.
This and \eqref{H-K0line} yield $H\cdot K_X\leq H^2=d$
and the rest of the argument is as in the first case.
\qed

We now bring in an additional hypothesis:
\begin{equation}
  \label{h:no-irp}
  X \mbox{ {\it is an irregular surface of general type
    with no irrational pencil}}.
\end{equation}

This hypothesis combined with the assumption $m_X\leq 5$ yield the
following charac\-terization.

\begin{lem}
  \label{l:4-5}
Let $m_X\leq 5$ and assume $X$ to be subject to \eqref{h:no-irp}.
Then $m_X=4$ or $5$.  Furthermore,

$\bullet$ if $m_X=4$, then $q(X)=2$ and 
\begin{equation}
  \label{K1}
  K_X = \mu H + F + \sum_{C\cdot H=2}m_C C + \sum_{l\cdot H=1} m_l l
\end{equation}
with $\mu\geq 1$, $F$ an effective divisor, and the sums run over
collections of irreducible conics ($C\cdot H=2$) and lines 
($l\cdot H=1$) contained in $X$;

$\bullet$ if $m_X=5$ and $X$ is minimal, then $q(X)\leq 3$.  Moreover,
if equality holds, then one
of the following can occur: 
\begin{enumerate}
\item[(i)]
  $\Omega_X$ is generated by global sections, and
  $K^2_X=6\chi(\Oo_X)$.
\item[(ii)]
  The canonical divisor has the following form
\[
  K_X = \mu H + F + \sum_{C\cdot H=2}m_C C + \sum_{l\cdot H=1} m_l l,
\]
where $\mu\geq 1$, $F$ is an effective divisor, and the sums are
as in \eqref{K1}.
\end{enumerate}
\end{lem}

\begin{rem*}
It is to be noticed that even if the coefficient $\mu$ above satisfies the
same condition in both cases, these conditions occur for different
geometric reasons.  The proof will show that if $m_X=4$, then
$K_X=E$ and $\mu$ correspnds to the sum, $\mu_{II}$, of the coefficients of
\typeII\ componenets of $E$, whereas if $m_X=5$ then $L=H$, \ie
$K_X=H+E$, and $\mu=\mu_{II}+1$.
\end{rem*}

\proof
The proof will be done in several steps.  We begin by showing in {\it
Step 1} that $m_X=4$ or $5$.  Then, we go on with the detailed study
of each value in {\it Step 2} and {\it Step 3}.  

\step Note that the hypothesis \eqref{h:no-irp} implies that the
irregularity $q(X)\geq 2$.  This and the fact that $X$ is of general
type yield $p_g(X)\geq 2$ as well.  Then, according to Lemma
\ref{l:EP-coh}, the coboundary homomorphism $\delta_X$ in
\eqref{coh-seq} is injective unless $m_X =5$, $p_g(X)\leq2$, and the
degree of $X$ in $\PP^4$ is at most $16$.  Hence, in this situation,
$p_g(X)=q(X)=2$.  Thus
we may assume that $\delta_X$ is injective, hence that
$H^1(\Theta_X(-K_X -(5- m_X)H))\neq 0$.

Let $\eta$ be a nonzero section of $\Nn(-K_X-(5-m_X)H)$ and let
$\xi=\delta_X(\eta)$ be its image in $H^1(\Theta_X(-K_X-(5-m_X)H))$.
Viewing it as an extension class in
$\Ext^1(\Omega_X,\Oo_X(-K_X-(5-m_X)H))$ via the natural identification
\[
  H^1(\Theta_X(-K_X-(5-m_X)H)) \cong 
  \Ext^1(\Omega_X,\Oo_X(-K_X-(5-m_X)H)),
\]
we obtain the following exact sequence of locally free sheaves on $X$,
\begin{equation}
  \label{ext1}
  0 \lra \Oo_X(-K_X-(5-m_X)H) \lra \Tt_\xi \lra \Omega_X 
  \lra 0.
\end{equation}
Observe that the divisor $(K_X+(5-m_X)H)$ is nef and big for $m_X\leq 5$.
This implies that $H^0(\Tt_\xi)=H^0(\Omega_X)$ and we
consider the smallest
saturated subsheaf $\Ff$ of $\Tt_\xi$ containing the image of the
evaluation map
\[
  H^0(\Tt_\xi)\otimes\Oo_X \lra \Tt_\xi.
\]
The hypothesis \eqref{h:no-irp} implies that the rank of $\Ff$ is at
least $2$.  On the other hand, from the defining sequence
\eqref{ext1}, it follows that $\det(\Tt_\xi)=\Oo_X(-(5-m_X)H)$. This
implies that $\Ff$ has rank $2$ unless $m_X =5$.  In this case, with
the additional assumption of $X$ being minimal, we still have $\Ff$ of
rank $2$, unless $\Ff = \Tt_\xi =\oplus^3\Oo_X$ (these equalities
follow immediately from $\det(\Tt_\xi) = \Oo_X$). This yields the case
{\it (i)} of the lemma.

Thus we may assume that $\Ff$ has rank $2$ and we arrive at a
diagram analogous to the one in \eqref{d:unstable},
\begin{equation}
  \label{d1:unstable}
  \UseTips
  \xymatrix@W=5mm@H=5mm{
    && 0 \ar[d] \\
    && \Ff \ar[d]\ar[rd]^{\phi_\xi} \\
    0 \ar[r] & \Oo_X(-K_X-(5-m_X)H) \ar[r] \ar[rd]_{\psi_\xi}&
    \Tt_\xi \ar[r]\ar[d] & \Omega_X \ar[r] & 0 \\
    && \Jj_Z(-L -(5-m_X)H) \ar[d] \\
    && 0
  }
\end{equation}
where $\Oo_X(L)=\det(\Ff)$.  By construction, the morphism $\phi_\xi$
is generically an isomorphism.  Arguing as in the proof of Lemma
\ref{l:decomp}, we obtain the decomposition of the canonical
divisor
\begin{equation}
  \label{K=LE}
  K_X = L + E
\end{equation}
where $L$ is given by the determinant of $\Ff$ and
$E=c_1(\coker(\phi_\xi))$ is a non-zero effective divisor.  We also
observe that the hypothesis \eqref{h:no-irp} implies that $L$ is
effective as well.

Let $e$ be a section of $\Oo_X(E)$ corresponding to the divisor $E$.
Similarly to the statement of Lemma~\ref{l:E-coh}, we have that $e$
annihilates the cohomology class $\xi$, \ie
\begin{equation}
  \label{e-xi}
  e\cdot\xi = 0 \textq{in}
  H^1(\Theta_X(E-K_X-(5-m_X)H)).
\end{equation}
From the commutative diagram
\begin{equation}
  \label{d:e-sq}
  \xymatrix{
    & H^0(\Nn_X(-D)) \ar[d]_e\ar[r]^{\delta_X} &
    H^1(\Theta_X(-D)) \ar[d]^e\\
    H^0(\Theta_{\PP^4}\otimes\Oo_X(E-D)) \ar[r] &
    H^0(\Nn_X(E-D)) \ar[r]^{\delta_X(E)} &
    H^1(\Theta_X(E-D))
  }
\end{equation}
where $D=K_X+(5-m_X)H$, it follows that
\[
  \delta_X(E)(e\cdot\eta)=e\cdot\delta_X(\eta)=e\cdot\xi=0,
\]
where the last equality is \eqref{e-xi}.  From this and the diagram
\eqref{d:e-sq}, it follows that the non-zero section $e\cdot\eta$ of
$\Nn_X(E-D)$ comes from
$H^0(\Theta_{\PP^4}\otimes\Oo_X(E-D))$.  Thus, using the
decomposition \eqref{K=LE}, we obtain
\[
  H^0(\Theta_{\PP^4}\otimes\Oo_X(-L-(5-m_X)H)) =
  H^0(\Theta_{\PP^4}\otimes\Oo_X(E-K_X-(5-m_X)H)) \neq 0.
\]

From the Euler sequence of $\Theta_{\PP^4}$, the non-vanishing of
$H^0(\Theta_{\PP^4}\otimes\Oo_X(-L-(5-m_X)H))$ may occur for one of
the following reasons: either
\begin{equation}
  \label{r1}
  H^0(\Oo_X(-L -(4-m_X)H)) \neq 0,
\end{equation}
or
\begin{equation}
  \label{r2}
  \ker\Big(
  H^1(\Oo_X(-L-(5-m_X)H)) \to
  H^0(\Oo_X(H))^\ast\otimes H^1(\Oo_X(-L-(4-m_X)H))
  \Big) \neq 0.
\end{equation}
It is obvious that \eqref{r1} implies that $m_X\geq 4$.  We claim that
the same condition is necessary for \eqref{r2}.  Indeed, let $h$ be a
nonzero section of $\Oo_X(H)$ and let $C_h$ be the corresponding
divisor in the linear system $|H|$.  Then, from the exact
sequence 
\begin{multline*}
  0 \lra \Oo_X(-L-(5-m_X)H) \lra \Oo_X(-L-(4-m_X)H) \lra\\
  \lra \Oo_{C_h}( -L-(4-m_X)H) \lra 0,
\end{multline*}
it follows that \eqref{r2} forces the non-vanishing of the cohomology
group $H^0(\Oo_{C_h}(-L-(4-m_X)H))$.  For this to hold, one must again
have $m_X\geq 4$.  This completes the proof of the first assertion of
the lemma.  In the rest of the proof we consider separately the two
possible values of $m_X$.

\step 
If $m_X=4$, then the above argument shows that for \eqref{r1}
(respectively for \eqref{r2}) to hold, $L$ must be $0$.  This together
with the definition of the sheaf $\Ff$ in \eqref{d1:unstable} imply
\[
\Ff = \oplus^2 \Oo_X.
\]
In particular, $q(X)=2$.

Next we turn to the formula for $K_X$ in \eqref{K1}.  For this, we use
the decomposition of $K_X$ given in \eqref{K=LE}.  Due to $L=0$, this
simply reads
\[
  K_X = E = E_I + E_{II},
\]
where the last equality is as in \eqref{ty-d}.  Furthermore from the
diagram \eqref{d1:unstable} with $m_X=4$, the conditions
characterizing the curves of \typeI\ and \typeII\ become
\begin{alignat}{4}
  &\text{{\bf\typeI:}} &\quad\label{typeI-a}
  &\deg Z_C \leq \deg\eta^\ast_C Z_C = 2 -2g(\Ctilde)-H\cdot C
  &&\hspace{4cm}\\
  &\text{{\bf\typeII:}} &\quad\label{typeII-a}
  &H^0(\Oo_C(C-H)) \neq 0
  &&\hspace{4cm}
\end{alignat}
where $\eta_C:\Ctilde\to C$, in \eqref{typeI-a}, denotes the
normalization of $C$.

From \eqref{typeI-a}, it follows that the irreducible components of
$E_I$ are either lines, $C\cdot H=1$, or conics, $C\cdot H=2$. 
From \eqref{typeII-a}, it follows that $C-H =F_C$ is an effective
divisor.  Hence, the canonical 
divisor $K_X$ satisfies   
\begin{equation}
  \label{formulaE}
  K_X = E = \mu_{II} H +
  \sum_{\substack{C\text{ of}\\\text{\typeII}}} m_C F_C
  + \sum_{C\cdot H=2} m_C C + \sum_{l\cdot H=1} m_l l,
\end{equation}
where $m_D$ stands for the multiplicity of a component $D$ in $K_X$
and
\[
  \mu = \mu_{II} = \sum_{C\text{ of \typeII}} m_C.
\]
This formula contains all the ingredients of \eqref{K1} except the
assertion about $\mu_{II}$. 
To see that, let us assume that 
$\mu_{II}=0$.  Then
\[
  K_X = \sum_{C\cdot H=2} m_C C + \sum_{l\cdot H=1} m_l l,
\]
\ie the canonical divisor is composed entirely of rational
curves.  Thus it is contracted to points by the Albanese map and
hence, the Zariski decomposition of $K_X$ has no positive part,
contradicting the fact that $X$ is of general type.

\step If $m_X=5$ and $X$ is minimal, then the construction of the
subsheaf $\Ff$ in \eqref{d1:unstable} goes through and we may assume
it to be of rank $2$ (otherwise the argument in {\it Step 1} shows
that {\it(i)} of the lemma holds).  This gives the decomposition of
the canonical divisor $K_X$ as in \eqref{K=LE}.  Continuing the
argument as in {\it Step 1}, we arrive at
\begin{equation}
  \label{eq:non-vanishing}
  H^0(\Theta_{\PP^4}\otimes\Oo_X(-L)) \neq 0.
\end{equation}

At this point we assume that the irregularity $q(X)\geq 3$ (otherwise
there is nothing to prove).  Consider the homomorphism
\begin{equation}
  \label{wedgeFtoL}
  \wedge^2 H^0(\Ff) \lra H^0(\Oo_X(L)).
\end{equation}
The hypothesis \eqref{h:no-irp} and the argument in the lemma of
Castelnuovo-de Franchis, see \cite[Proposition X.9]{Be},
imply
\begin{equation}
  \label{dimL}
  h^0(\Oo_X(L)) \geq 2q(X) -3 \geq 3.
\end{equation}
We claim that this leads to
\begin{equation}
  \label{L=H}
  \Oo_X(L) = \Oo_X(H).
\end{equation}
Let us assume this and complete our argument.  The inequality
\eqref{dimL} combined with \eqref{L=H} yield
\begin{equation}
  \label{b1}
  2q(X) -3 \leq h^0(\Oo_X(L)) = h^0(\Oo_X(H)) = 5,
\end{equation}
\ie $q(X)\leq 4$.  We will now rule out the case $q(X)=4$.  Indeed, if
$q(X)=4$ holds, the homomorphism in \eqref{wedgeFtoL} has the form
\begin{equation}
  \label{Plucker}
  \wedge^2 H^0(\Ff) \lra H^0(\Oo_X(H))
\end{equation}
and must be surjective with a $1$-dimensional kernel.  Hence $\Ff$ is
generated by global sections and it maps $X$ into the Grassmannian
$\Gr(1,3)=\Gr(1,\PP(H^0(\Ff)^\ast))$ of lines in
$\PP(H^0(\Ff)^\ast)=\PP^3$.  Composing this map with the Pl\"ucker map
of the Grassmannian, we get
\[
  X \longrightarrow \Gr(1,3) \hookrightarrow
  \PP(\wedge^2 H^0(\Ff)^\ast) = \PP^5
\] 
Furthermore $X\subset\PP^4$ is contained in the intersection of the
Pl\"ucker embedding of the Grassmannian $ Gr(1,3)$ with the hyperplane
in $\PP(\wedge^2 H^0(\Ff)^\ast)$ corresponding to the kernel of the
homomorphism in \eqref{Plucker}.  The image of $\Gr(1,3)$ under the
Pl\"ucker map is a quadric.  Hence $X\subset\PP^4$ is contained in a
quadric which is impossible by the first assertion of the lemma.

Thus $q(X)=3$ and the decomposition in \eqref{K=LE} reads as follows:
\begin{equation}
  \label{K=HE}
  K_X = H + E = H + E_{II} + E_{I}
\end{equation}
where we use again the decomposition of $E$ into types.  Furthermore,
the irreducible components of $E$ are subject to the same conditions
as in \eqref{typeI-a} and \eqref{typeII-a}.  This leads to the formula
for $E$ as in \eqref{formulaE}.  Substituting it into \eqref{K=HE}
yields
\[
K_X = (\mu_{II}+1)H +
\sum_{\substack{C\text{ of}\\\text{\typeII}}} m_C F_C +
\sum_{C\cdot H=2} m_C C + \sum_{l\cdot H=1} m_l l
\]
where the notation have the same meaning as in \eqref{formulaE}.
\qed

\likeproof[Proof of \eqref{L=H}]
For $m_X =5$ the conditions \eqref{r1} and \eqref{r2}  
become respectively $H^0(\Oo_X(H-L))\neq 0$ and 
\begin{equation}
  \label{R2}
  \ker\Big(
  H^1(\Oo_X(-L)) \lra H^0(\Oo_X(H))^\ast\otimes H^1(\Oo_X(H-L))
  \Big) \neq 0.
\end{equation}

\paragraph{Claim}
$H^0(\Oo_X(H-L))\neq 0$.

\bigskip

This claim implies the equality \eqref{L=H}. 
Indeed, assume  
$H-L=\Gamma\neq 0$. Then we have
\[
  h^0(\Oo_X(H-\Gamma)) = h^0(\Oo_X(L)) \geq 3,
\]
where the last inequality comes from \eqref{dimL}.  This implies that
the equality must hold and that $\Gamma$ is a line in $\PP^4$.  Thus
\begin{equation}
  \label{formH}
  H = L + \Gamma = K-E +\Gamma
\end{equation}
Taking the intersection with $\Gamma$, we have
\[
  1 = \Gamma\cdot H = 
  \Gamma\cdot(K_X+\Gamma) - \Gamma\cdot E = -2-\Gamma\cdot E.
\]
From this it follows that $E\cdot\Gamma=-3$ which contradicts
Remark~\ref{Knonample}.

\bigskip

We now turn to the proof of the claim, $H^0(\Oo_X(H-L))\neq 0$.  From
\eqref{dimL}, it follows that
\begin{equation}
  \label{L=M+F}
  L=M+F
\end{equation}
where $M$ (resp. $F$) is the moving (resp. fixed) part of $L$.  In
particular, $M$ is nef, and since $X$ has no irrational pencil, it is
also big.  The non-vanishing of the group
$H^0(\Theta_{\PP^4}\otimes\Oo_X(-L))$ (see \eqref{eq:non-vanishing})
implies that
\[
  H^0(\Theta_{\PP^4}\otimes\Oo_X(-M))\neq0 
\]
as well.  From the Euler sequence of $\Theta_{\PP^4}$ it follows that
this group is the middle term of the exact sequence
\begin{equation*}
  H^0(\Oo_X(H))^\ast\otimes H^0(\Oo_X(H-M)) \lra
  H^0(\Theta_{\PP^4}\otimes\Oo_X(-M)) \lra H^1(\Oo_X(-M)).
\end{equation*}
The fact that $M$ is nef and big yields $H^1(\Oo_X(-M))=0$.  Hence
$H^0(\Oo_X(H-M))\neq 0$.  Using this non-vanishing and the decomposition
\eqref{L=M+F}, $L=M+F$, we see that 
\begin{equation}
  \label{eq:Ldecomp}
  L = H+D,
\end{equation}
with $D$ an effective divisor.  Indeed, $H^0(\Oo_X(H-M))\neq 0$ if
either $M=H$, and $D=F$ in \eqref{eq:Ldecomp}, or $H-M=l$, with $l$ a
line since $h^0(\Oo_X(H-l))=h^0(\Oo_X(M))=h^0(\Oo_X(L))\geq 3$.  Now
$l\cdot E\geq -2$ (see Remark \ref{Knonample}), hence
\[
  0 \geq
  l\cdot (K_X+l-E) = l\cdot (L+l) =
  l\cdot (H+L-M) = l\cdot (H+F) = 1+l\cdot F.
\]
It follows that $l$ is contained in $F$ and we take $D=F-l$ to obtain
the decomposition in \eqref{eq:Ldecomp}.

We now use \eqref{eq:Ldecomp} to show that the kernel in \eqref{r2}
vanishes.  Indeed, a non-trivial element of this kernel gives rise to
the non-trivial kernel for the cup product
\[
  H^1(\Oo_X(-L)) \stackrel{h}{\lra} H^1(\Oo_X(H-L))
\]
for every $h\in H^0(\Oo_X(H))$.
Thus $H^0(\Oo_{C_h}(H-L))\neq 0$, where $C_h=\{h=0\}$ is the divisor
corresponding to $h$.  Substituting for $L$ the expression
\eqref{eq:Ldecomp}, we deduce that
$H^0(\Oo_{C}(-D))\neq 0$  for
every divisor $C$ in the linear system $|H|$.  Since $D$ is effective,
we obtain that $D=0$ and thus $L=H$.  This implies the vanishing of
$H^1(\Oo_{X}(-L))$ in \eqref{R2}, hence $H^0(\Oo_X(H-L))\neq 0$.
\qed

\likeproof[Proof of Theorem \ref{bound-q}]
Let $d(\chi,\beta)$ be as in Lemma \ref{l:EP} and let
$X\in\sS_d^\circ(\chi,\beta)$.  By Lemma~\ref{l:EP}, $X$ is contained
in a $3$-fold of degree $m_X\leq 5$.  From Lemma~\ref{l:4-5} it
follows that $m_X$ is either $4$ or $5$, and that $q(X)\leq 3$.
Furthermore, if $q(X)=3$, then {\it (i)} or {\it (ii)} of Lemma
\ref{l:4-5} must hold.
\qed

\begin{rem*}
The well-known example of Horrocks and Mumford in \cite{HM} is
essentially the only known surface in $\PP^4$ of irregularity $2$.
Though Theorem~\ref{bound-q} does not rule out the possibility of
surfaces in $\PP^4$ with irregularity $3$, it shows, following
\cite[Proposition 4.1]{De}, that such a hypothetical
surface\footnote{The condition {\it (ii)} of Lemma \ref{l:4-5} could
be envisaged as a degenerate case of {\it (i)}.}  could be a divisor
in an Abelian variety of dimension $3$.
\end{rem*}

Using similar reasoning, we can restrict the topology of surfaces of
high degree and bounded holomorphic Euler characteristic in $\PP^4$.

\begin{thm}
  \label{t:topCondition}
Given a positive integer $n$, every surface $X$ in $\sS_d(\chi,\beta)$
has negative topological index, \ie $\alpha_X>1/2$, provided
$d>d(\chi,\beta)$.  (The notation $\sS_d(\chi,\beta)$ and
$d(\chi,\beta)$ are to be found in \eqref{eq:Sdn} and Lemma~\ref{l:EP}
respectively.)
\end{thm}

\proof
By Lemma~\ref{l:EP}, $X$ is contained in a $3$-fold of degree 
$m_X\leq 5$.  Following the arguments in Lemma~\ref{l:EP-coh} and
Lemma~\ref{l:4-5}, we deduce that $H^1(\Theta_X(-K_X))\neq0$.  This
implies the decomposition of $K_X=L+E$ and that
$H^0(\Theta_{\PP^4}\otimes\Oo_X(-L))\neq0$.  The latter condition
yields
\begin{equation}
  \label{eq:HL}
  H\cdot L \leq H^2 = d.
\end{equation}

\paragraph{Claim}
If $\alpha_X\leq1/2$, then $E\cdot H\leq L\cdot H$.

\medskip

To justify the claim, we assume the opposite, $E\cdot H>L\cdot H$.
From this and Lemma~\ref{l:positiveCone}, 1), we deduce 
\[
  (E-L)\cdot (K_X+\lambda H) = 0
\]
for some $\lambda>0$.  By the Hodge index and
Lemma~\ref{l:positiveCone}, 2), it follows that $E$ is numerically
equivalent to $L$, which contradicts Lemma~\ref{l:positiveCone}, 1).

\medskip

To end the proof of the theorem, we combine the claim and
\eqref{eq:HL} to deduce
\[
H\cdot K_X = H\cdot (L+E) \leq 2H\cdot L \leq 2d.
\]
This inequality and \eqref{eq:d-and-Chern} yield
\[
2d \geq H\cdot K_X = \frac{1}{5}(d^2 -10d-2(\beta-6)\chi).
\]
Hence, when the corresponding discriminant is not negative, 
\[
d \leq 10+\sqrt{100+2(\beta-6)\chi} < d(\chi,\beta).
\]
\qed

\bigskip

\begin{tabular}[r]{p{22ex}rp{3ex}r}\small
  &  Daniel Naie && Igor Reider \\
  & Universit\'e d'Angers && Universit\'e d'Angers \\
  & 2, Bd Lavoisier && 2, Bd Lavoisier \\
  & France && France \\
  &daniel.naie@univ-angers.fr && igor.reider@univ-angers.fr
\end{tabular}


\begin{thebibliography}{99}
\bibitem{B-C}
  {\sc E.~Ballico, L.~Chiantini},
  On smooth subcanonical varieties of codimension 2 in $\PP^n$, $n\geq4$.
  {\it Ann. Mat. Pura Appl. 135\/}  (1984), 99–-117. 
\bibitem{Be}
  {\sc A.~Beauville},
  {\it Surfaces alg\'ebriques complexes}, 
  Ast\'erisque, No. 54. Soci\'et\'e Math\'ematique de France, Paris, 1978.
\bibitem{Bog}
  {\sc  F.~A.~Bogomolov},
  Unstable vector bundles and curves on surfaces.
  {\it Proceedings of the International Congress of Mathematicians}
 (Helsinki, 1978), 517--524.
\bibitem{De}
  {\sc O. Debarre},
  Th\'eor\`emes de connexit\'e et vari\'et\'es
  ab\'eliennes. 
  {\it Amer. J. Math. 117\/} (1995), no. 3, 787--805.
\bibitem{DeSch}
  {\sc W.~Decker, F.-O.~ Schreyer},
  Non-general type surfaces in $\PP^4$: some remarks on bounds and
  constructions. 
  {\it J. Symbolic Comput. 29\/} (2000), no. 4-5, 545--582. 
\bibitem{E-P}
  {\sc G.~Ellingsrud, C.~Peskine},
  Sur les surfaces lisses de $\PP^4$, 
  {\it Invent.math.~ 95\/} (1989), 1--11.
\bibitem{G-P}
  {\sc  L.~Gruson, C.~Peskine},
  Genre des courbes de l'espace projectif.
  {\it Algebraic Geometry, LNM 687\/}, Springer 1977, 31--59.
\bibitem{H}
  {\sc R.~Hartshorne}, 
  {\it Algebraic Geometry}, Springer, 1977.
\bibitem{HM}
  {\sc G.~Horrocks, D.~Mumford},
  A rank $2$ vector bundle on $\PP^{4}$ with $15,000$ symmetries.
  {\it Topology 12\/} (1973), 63--81. 
\bibitem{KS}
  {\sc K.~Kodaira, D.~C.~Spencer},
  On deformations of complex analytic structures, I-II.
  {\it Annals of Math.~ 67\/} (1958) 328--466. 
\bibitem{L}
  {\sc J. Lipman},
  Free derivation modules on algebraic varieties.
  {\it American Journal of Math. 87\/}, (1965), 874--898.
\bibitem{Miy}
  {\sc Y.~Miyaoka},
  The maximal number of quotient singularities on surfaces with given
  numerical invariants.
  {\it Math. Ann.~ 268\/} (1984), 159--171.
\bibitem{OSS}
  {\sc Ch.~Okonek, M.~Schneider, H.~Spindler}, 
  {\it Vector bundles on complex projective spaces},
  Progress in Mathematics, 3. Birkh\"auser, Boston, Mass., 1980.
\bibitem{R}
  {\sc I. Reider},
  Geography and the number of moduli of surfaces of general type.
  {\it Asian J.Math. 9\/} (2005), 407--448.
\bibitem{Ro}
  {\sc L.~Roth}, 
  On the projective classification of surfaces,
  {\it Proc. London Math. Soc. 42\/} (1937), 142--170.
\end{thebibliography}
\end{document}